\documentclass[12pt]{article}

\usepackage{amsfonts}
\usepackage{graphicx}
\usepackage{enumerate}
\usepackage[shortlabels]{enumitem}
\usepackage{mathtools}
\usepackage{hyperref}
\usepackage{color}
\usepackage{xcolor}
\usepackage{fullpage}
\usepackage{tikz}
\usepackage{cite}  

\usepackage{amsmath, amssymb,latexsym}
\usepackage{soul}
\usepackage{float}
\usepackage{comment}

\def\C{\mathbb{C}}
\def\F{\mathbb{F}}

\def\N{\mathbb{N}}

\def\Z{\mathbb{Z}}

\DeclareMathOperator{\diag}{diag}

\DeclareMathOperator{\rank}{rank}

\newcommand{\fpartial}{\hat{\partial}}

\newcommand{\ideal}[1]{\langle #1 \rangle}

\newtheorem{theorem}{Theorem}[section]
\newtheorem{proposition}[theorem]{Proposition}

\newtheorem{lemma}[theorem]{Lemma}
\newtheorem{definition}[theorem]{Definition}
\newtheorem{corollary}[theorem]{Corollary}
\newtheorem{remark}[theorem]{{\sc Remark}}
\newtheorem{example}[theorem]{Example}
\newtheorem{problem}[theorem]{Problem}

\usepackage[normalem]{ulem}
\definecolor{brilliantrose}{rgb}{1.0, 0.33, 0.64}
\definecolor{myviolet}{rgb}{0.21, 0.0, 0.85}
\definecolor{amethyst}{rgb}{0.6, 0.4, 0.8}
\definecolor{carrotorange}{rgb}{0.93, 0.57, 0.13}
\definecolor{cutepink}{rgb}{1.0, 0.2, 0.6}

\title{The Jordan canonical form of\\ the Fr\'{e}chet derivative of a matrix function\\
and the bivariate Jordan problem}
\author{Vanni Noferini\thanks{Aalto University, Department of Mathematics and Systems Analysis, P.O. Box 11100, FI-00076, Aalto, Finland. Supported by a Research Council of Finland grant (decision number 370932). Email: vanni.noferini@aalto.fi}}
\date{6 May 2026}

\begin{document}

\maketitle

\begin{abstract}
 Let $\F$ be an algebraically closed field of characteristic $0$. Given a square matrix $A \in \F^{n \times n}$ and a polynomial $f \in \F[w]$, we determine the Jordan canonical form of the formal Fr\'{e}chet derivative of $f(A)$, in terms of that of $A$ and of $f$. When $\F\subseteq \C$, via Hermite interpolation, our result provides a solution to [N.J. Higham, \emph{Functions of Matrices: Theory and Computation}, Research Problem 3.11]. A generalization consists of finding the Jordan canonical form of linear combinations of Kronecker products of powers of two square matrices, i.e., $\sum_{i,j} a_{ij} (X^i \otimes Y^j)$. For this generalization, we provide some new partial results, including a partial solution under certain assumptions and general bounds on the number and the sizes of Jordan blocks.
\end{abstract}


\bigskip


\section{Motivation and statement of the problem}

Let $(U,\|\cdot\|_U)$, $(V,\|\cdot\|_V)$ be Banach spaces. In functional analysis \cite[Chapter 2]{Berger_77}, given an open set $\Omega \subset V$, a function $f: U \rightarrow V$ is said to be Fr\'{e}chet differentiable at $x \in \Omega$ if a continuous linear operator $Df_{w}$ can be found such that

\[ \lim_{\|e\|_U \rightarrow 0 } \frac{\|f(w+e)-f(w)-Df_w(e)\|_V}{\|e\|_U} =0. \]

The operator $Df_w$, when it exists, is unique and is known as the Fr\'{e}chet derivative of $f$ at $w$. When $U=V=\mathbb{C}^{n \times n}$, a case of particular interest is when $f$ is a matrix function  in the linear algebraic sense, which is a more restrictive notion than just applying the usual definition of a function to a matrix \cite{Higham_book, HJ}. Since $U$ and $V$ are now normed vector spaces of finite dimension $n^2$ over $\mathbb{C}$, it is possible to represent the linear map $E \rightarrow Df_W(E)$ as an $n^2 \times n^2$ matrix. Clearly, the representation depends on the choice of a basis for $\mathbb{C}^{n \times n}$. The most common choice is induced by the canonical basis of $\mathbb{C}^n$ via the two formal operations of Kronecker product and vec operator \cite[p. 331]{Higham_book}. Such representation of $Df_W$ is sometimes referred to as the Kronecker form of the Fr\'{e}chet derivative. Furthermore, if $f$ is regular enough, say, analytic in an open set including the eigenvalues of $W$, then there exists a  polynomial $p$ (depending on both $f$ and $W$) such that $Df_W = Dp_W$ \cite[Theorem 6.6.14]{HJ}; see also the discussion in \cite[Section 3.2]{Higham_book}. Therefore studying the case of polynomial functions suffices to cover a vast class of functions. It can be shown \cite[Problem 3.6]{Higham_book} that for the polynomial function $\displaystyle f(W)=\sum_{i=-1}^{{L}} f_i W^{i+1}$
the Kronecker form of $Df_W$ is the matrix
\begin{equation}\label{eq:kron}
    K(f,W) = \sum_{i=0}^{{L}} f_i \sum_{j=0}^i (W^T)^j \otimes W^{i-j} ,
\end{equation} 
where $A \otimes B$ denotes the Kronecker product  of the two matrices $A$ and $B$.

This setting motivates the following problem. For $0 \leq i,j, \leq k$ let $a_{ij} \in \mathbb{F}$ where $\mathbb{F}$ is an algebraically closed field of  characteristic zero. Consider the bivariate polynomial 
\begin{equation}\label{pscal} p(x,y)= \sum_{0 \leq i,j \leq k} a_{ij} x^i y^j \in \F[x,y].\end{equation}

Given two square matrices $X \in \mathbb{F}^{m \times m}$, $Y \in \mathbb{F}^{n \times n}$, we study the matrix

\begin{equation}\label{pmat} 
P(X,Y) = \sum_{0 \leq i,j \leq k} a_{ij} (X^i \otimes Y^j)  \ \in \mathbb{F}^{mn \times mn}.
\end{equation}

\begin{remark}\label{rem:theycommute}
Fix a pair of square matrices $X \in \F^{m \times m},Y \in \F^{n \times n}$. The set of all matrices of the form \eqref{pmat}, equipped with matrix addition and matrix multiplication, is a commutative ring $\mathcal{K}_{X,Y} \subseteq \F^{mn \times mn}$. Moreover, if we denote by $\mu_M(t)$ the minimal polynomial of the square matrix $M$ {and by $\langle a(x,y),b(x,y) \rangle$ the ideal generated by $a(x,y),b(x,y) \in \F[x,y]$}, then by the first isomorphism theorem \[ \mathcal{K}_{X,Y} \cong \F[x,y] / \langle \mu_X(x),\mu_Y(y) \rangle,\]
via the ring isomorphism $p(x,y) \mapsto P(X,Y)$.
\end{remark}

We are especially interested in describing the Jordan canonical form (JCF) of \eqref{pmat} in terms of the scalar polynomial \eqref{pscal} and the JCFs of $X$ and $Y$. It is clear that $p(\lambda,\mu)$ is an eigenvalue of $\eqref{pmat}$ if and only if $\lambda$ is an eigenvalue of $X$ and $\mu$ is an eigenvalue of $Y$. Hence, the only nontrivial part of the computation of the JCF of \eqref{pmat} is to determine the sizes of each Jordan block.  We call this problem the \emph{bivariate Jordan problem}. Some specific instances of the bivariate Jordan problem have been solved long ago. If the only nonzero coefficient in \eqref{pmat} is $a_{11}=1$, we obtain the Kronecker product of two matrices, whose JCF has been thoroughly studied by various authors \cite{Aitken_34, Brualdi_85, Littlewood_35, Marcus_75, Roth_34}, and for this case the solution is well known. Another important special case is the Kronecker sum of two matrices: $a_{10}=a_{01}=1$, $a_{ij}=0$ otherwise. Kronecker sums are often used and a complete characterization of their JCF first appeared in \cite[Theorem 2]{Roth_34}.
 
 More generally, the bivariate Jordan problem was studied by C. Norman \cite{Norman_97}, who fully analyzed the  case $m,n \leq 3$ and gave partial results for arbitrary $m,n$. Our approach is mostly complementary to \cite{Norman_97}, although there is a partial superposition of the techniques and the results. Another related paper is \cite{Kressner_14}, where D. Kressner studied (for $\F=\C$) the class of linear operators $\displaystyle C \mapsto \sum_{0\leq i,j \leq k}a_{ij} Y^j C (X^T)^i$; via the vec operator, these are precisely the linear functions on $\F^{n \times m}$ whose matrix representation is \eqref{pmat}. While \cite{Kressner_14} discusses spectral properties, it touches upon the JCF only marginally, and it explicitly refers to  \cite{Norman_97} for more details. Some of our results can also be seen as extensions of the theory in \cite{Kressner_14}.
 
 Motivated by \cite{Higham_book}, we are especially interested in the case where $P(X,Y)$ mimics the form \eqref{eq:kron}. Namely, for some $f_i \in \F$, we are interested in the special case
\begin{equation}\label{eq:frechet}
   P(X,Y) = \sum_{i=0}^{{L}} f_i \left(\sum_{j=0}^i X^j \otimes Y^{i-j}\right) \ \in \F^{mn \times mn}.
\end{equation}
When $X=Y^T$, it is clear that \eqref{eq:frechet} generalizes \eqref{eq:kron} from $\C$ to a more general field $\F$, and in that case we call it the \emph{formal Fr\'{e}chet derivative} at $Y$ of the polynomial function \[ f(W) : \F^{n \times n} \rightarrow \F^{n \times n}, \qquad W \mapsto f(W)=\sum_{i=-1}^{{L}} f_i W^{i+1};\] we also call \emph{Fr\'{e}chet-Jordan problem} the special case of the bivariate Jordan problem where \eqref{pmat} takes the form \eqref{eq:frechet}. Our motivation is to answer \cite[Research Problem $3.11$]{Higham_book}, in which N. Higham asked this question for the Kronecker form of the Fr\'{e}chet derivative of a matrix function. In other words, \cite[Research Problem $3.11$]{Higham_book} is the Fr\'{e}chet-Jordan problem for $\mathbb{F}=\mathbb{C}$ and $X=Y^T$. To our knowledge, \cite[Research Problem $3.11$]{Higham_book} was until now still an open question. The present paper answers it.

As already noted in \cite{Norman_97}, there is a technical assumption about the triple $(p,X,Y)$ that can simplify the bivariate Jordan problem. We state such assumption locally for each pair $(\lambda,\mu)$ of eigenvalues of $X$ and $Y$, as inequality \eqref{eq:nonzerograd} in Subsection \ref{sec:generic}.  If we fix $m,n$, and an upper bound on the total degree of \eqref{pscal}, so that the problem's parameters lie in $\F^N$ for some $N$, then \eqref{eq:nonzerograd} is a generic assumption in the Zariski topology of $\F^N$. On the other hand, the degenerate case when \eqref{eq:nonzerograd} fails is a wild combinatorial problem, whose complete solution appears difficult. The main contributions of this paper are: 
\begin{itemize}
    \item[(i)] For an arbitrary bivariate polynomial \eqref{pscal}, and for each pair of Jordan blocks in the JCF of $X,Y$ for which \eqref{eq:nonzerograd} holds, we fully characterize the associated Jordan blocks in the JCF of \eqref{pmat}; in particular, this solves the bivariate Jordan problem when \eqref{eq:nonzerograd} holds for every eigenvalue pair of $X$ and $Y$. For the eigenvalue pairs for which \eqref{eq:nonzerograd} fails, we  instead provide bounds on the number and the sizes of the associated Jordan blocks in the JCF of \eqref{pmat}.
    \item[(ii)] For the special case \eqref{eq:frechet}, we give a full general solution to the Fr\'{e}chet-Jordan problem, regardless of whether \eqref{eq:nonzerograd} holds or fails. Our solution is by a closed formula when $X$ and $Y$ do not share eigenvalues, whereas for shared eigenvalues it reduces the problem to computing the ranks of certain banded integer Toeplitz matrices (of much smaller size than \eqref{eq:frechet}). By taking $\F=\C$ and $X=Y^T$ in \eqref{eq:frechet}, this in particular solves \cite[Research Problem 3.11]{Higham_book}.
\end{itemize}

We structure our exposition as follows. We make some preliminary analysis of the bivariate Jordan problem in Section \ref{sec:simp}. In Subsection \ref{sec:generic} we solve the bivariate Jordan problem assuming that \eqref{eq:nonzerograd} holds for all eigenvalue pairs ($\lambda,\mu)$, which as discussed happens for a non-empty and Zariski-open subset of problem instances of bounded size and degree. In Subsection \ref{sec:bounds} we give general bounds for the degenerate case when the assumption is not true. In Section \ref{sec:further}, we fully solve the Fr\'{e}chet-Jordan problem in the sense mentioned above. Finally, some conclusions are drawn in Section \ref{sec:conclusions}.

\section{Simplifying preliminaries}\label{sec:simp}

Suppose that $X=S J_X S^{-1}$ and $Y = R J_Y R^{-1}$ for certain nonsingular matrices $S$ and $R$, with $J_X, J_Y$ both in JCF. Then by \cite[Lemma 4.2.10]{HJ} $(S^{-1} \otimes R^{-1}) (X^i \otimes Y^j) (S \otimes R) = J_X ^i \otimes J_Y^j$. Therefore, $P(X,Y)$ and $P(J_X,J_Y)$ are similar, which means that there is no loss of generality in taking both $X$ and $Y$ in JCF. Furthermore, it suffices to separately analyze each pair $(X_i,Y_j)$ where $X_i$ (resp. $Y_j$) is a Jordan block of $X$ (resp. $Y$).

\begin{lemma}\cite[Corollary 4.3.16]{HJ}
Suppose that $X=\bigoplus_i X_i$ and $Y=\bigoplus_j Y_j$ are block diagonal. Then $X \otimes Y$ is similar to $\bigoplus_{i,j} (X_i \otimes Y_j)$. Moreover, the similarity matrix is a permutation matrix that only depends on the sizes of the blocks $X_i$ and $Y_j$.
\end{lemma}

A matrix $X$ in JCF is block-diagonal, with upper triangular Toeplitz blocks. Thus, $X^i$ is also block-diagonal with the same size of the (still upper triangular Toeplitz) blocks. Therefore, since $X$ and $Y$ are in Jordan form, $P(X,Y)$ is similar to the direct sum of blocks of the form $P(J_m(\lambda),J_n(\mu))$, where $J_m(\lambda)$ is any Jordan block in $X$ while $J_n(\mu)$ is any Jordan block in $Y$. Hence, it suffices to solve the bivariate Jordan problem when the matrix arguments of $P(\cdot,\cdot)$ are the single Jordan blocks
\[J_m(\lambda)= \lambda I_m + N_m, \qquad J_n(\mu)= \mu I_n + N_n. \]
Here and below, $N_k=J_k(0)$ denotes a nilpotent Jordan block of size $k$. Denoting formal partial derivatives with respect to $x$ by $\frac{\partial}{\partial x}$, it can be checked by direct computation that $J_m(\lambda)^i$ is the upper triangular Toeplitz matrix whose first row is {$\displaystyle \begin{bmatrix}
x^i & \frac{\partial x^i}{\partial x} & \frac{1}{2!}\frac{\partial^2 x^i}{\partial x^2} & \dots & \frac{1}{(m-1)!}\frac{\partial^{m-1} x^i}{\partial x^{m-1}}
\end{bmatrix}$} evaluated at $x=\lambda$; an analogous result holds for $J_n(\mu)^j$. By the definition of the Kronecker product one has that $J_m(\lambda)^i \otimes J_n(\mu)^j$ is equal to the evaluation at $x=\lambda$ of the matrix

\[ { \begin{bmatrix}
x^i J_n(\mu)^j & \frac{\partial x^i}{\partial x} J_n(\mu)^j & \frac{1}{2!}\frac{\partial^2 x^i}{\partial x^2} J_n(\mu)^j & \dots & \frac{1}{(m-1)!}\frac{\partial^{m-1} x^i}{\partial x^{m-1}} J_n(\mu)^j\\
0 & x^i J_n(\mu)^j & \frac{\partial x^i}{\partial x} J_n(\mu)^j & \dots & \frac{1}{(m-2)!}\frac{\partial^{m-2} x^i}{\partial x^{m-2}} J_n(\mu)^j\\
\vdots & \ddots & \ddots & \ddots & \vdots \\
\vdots & & \ddots & \ddots & \frac{\partial x^i}{\partial x} J_n(\mu)^j \\
 0 & \dots & \dots & 0 & x^i J_n(\mu)^j
\end{bmatrix}. } \]
By linearity we obtain a formula for $P(J_m(\lambda),J_n(\mu))$, which is a block-Toeplitz Toeplitz-block (bTTb) matrix that we describe below, omitting for simplicity the dependence on $(\lambda,\mu)$ of the evaluations of $p$ and its formal Hasse derivatives. These  are the $\F$-linear operators $\displaystyle \fpartial_{x^\alpha}  := \frac{1}{\alpha!}\frac{\partial^\alpha }{\partial x^\alpha}$ and $\displaystyle \fpartial_{y^\beta}  := \frac{1}{\beta!}\frac{\partial^\beta }{\partial y^\beta}$ such that
\[ \fpartial_{x^\alpha} \fpartial_{y^\beta} x^i y^j = \begin{cases}
   \binom{i}{\alpha} \binom{j}{\beta} x^{i -\alpha} y^{j-\beta} \ &\mathrm{if} \ \alpha \leq i \ \mathrm{and} \ \beta \leq j,\\
   0 \ &\mathrm{otherwise}.
\end{cases} \]

Within this formalism, we have

\begin{equation}\label{ourmatrix}
P(J_m(\lambda),J_n(\mu))=\begin{bmatrix}
B & \fpartial_x B & \dots & \dots & \fpartial_{x^{m-1}} B\\
0 & B & \fpartial_x B & \dots & \fpartial_{x^{m-2}} B\\
\vdots & \ddots & \ddots & \ddots & \vdots \\
\vdots & & \ddots & \ddots & \fpartial_x B \\
0 & 0 & \dots & 0 & B\end{bmatrix},
  \ \  B= \begin{bmatrix}
p & \fpartial_y p & \dots & \dots & \fpartial_{y^{n-1}} p\\
0 & p & \fpartial_y p & \dots & \fpartial_{y^{n-2}} p\\
\vdots & \ddots & \ddots & \ddots & \vdots\\
\vdots & & \ddots & \ddots & \fpartial_y p \\
0 & 0 & \dots & 0 & p
\end{bmatrix}. 
\end{equation} 

Lemma \ref{switch} below will let us swap the roles of the pairs $(x,m)$ and $(y,n)$ whenever this is convenient for our analysis. 

\begin{lemma}\label{switch}
Let $\displaystyle q(x,y):=p(y,x)$. Then $P(X,Y)$ is permutation similar to $Q(Y,X)$. The similarity matrix $\Pi$ can be written by columns as 
\[ \Pi = [e_1, e_{n+1},\dots, e_{(m-1)n+1}, e_2, e_{n+2}, \dots, e_{(m-1)n+2}, \dots, \dots, e_n, e_{2n}, \dots, e_{mn}],   \]
where $\{e_i\}_{i=1,\dots,mn}$ is the canonical basis of $\mathbb{F}^{mn}$.
\end{lemma}

\begin{proof}
Recall \cite[Corollary 4.3.10]{HJ} that, for any two square matrices $X$ and $Y$, $X \otimes Y$ and $Y \otimes X$ are permutation similar and the similarity matrix $\Pi$ only depends on the sizes of $X$ and $Y$. Comparing \eqref{pscal} with
$\displaystyle q(x,y) = \sum_{0 \leq i,j \leq k} a_{ji} x^i y^j$, we get that 
$\displaystyle Q(Y,X)= \sum_{0 \leq i,j \leq k} a_{ji} Y^i \otimes X^j = \sum_{0 \leq i,j \leq k} a_{ij} Y^j \otimes X^i = \Pi^T P(X,Y) \Pi$.
The explicit form of $\Pi$ is given in \cite[Equation 4.3.9b]{HJ}: It is straightforward to check that the formula given there is equivalent to the one given here.
\end{proof}

The problem is now to obtain the JCF of the bTTb matrix \eqref{ourmatrix}, or equivalently, as the only eigenvalue is obviously $p(\lambda,\mu)$, the ranks of the matrices $(P(J_m(\lambda),J_n(\mu))-p(\lambda,\mu)I_{mn})^k$, $k=1,\dots,\nu$, where $\nu$ {is} the smallest value of $k$ for which the rank is equal to $0$. Since the value of $p(\lambda,\mu)$ is irrelevant to the answer, we may also assume with no loss of generality that it is $0$, or in other words that $(P(J_m(\lambda),J_n(\mu))$ is nilpotent.

\begin{lemma}\label{lem:known}\cite[Section 0.9.7]{HJunder}
The set of $n \times n$ upper triangular Toeplitz matrices with elements in a commutative ring, equipped with matrix addition and matrix multiplication, is itself a commutative ring.
\end{lemma}

We now denote by $R$ the commutative ring of $n \times n$ upper triangular Toeplitz matrices with elements in $\mathbb{F}$ and by $\hat{R}$ the commutative ring of $m \times m$ upper triangular Toeplitz matrices with elements in $R$, that is, $\hat{R}$ is the ring of block upper triangular Toeplitz matrices with upper triangular Toeplitz blocks. The following observations are useful.
\begin{itemize}
\item The ring isomorphisms $R \cong \F[x]/\ideal{x^n}$ and $\hat{R} \cong R[y]/\ideal{y^m}\cong \F[x,y]/\ideal{x^n,y^m}$ hold.
\item $A \in R$ is invertible over $\F$ if and only if $A_{11} \neq 0$, while $Z \in \hat{R}$ is invertible over $\mathbb{F}$ if and only if $Z_{11} \in R$ is.
\item \cite{GS72} If $A \in R$ is invertible over $\F$ then it is a unit of $R$, i.e., $A^{-1} \in R$. Moreover, $Z \in \hat{R}$ is invertible over $\mathbb{F}$ if and only if it is invertible over $R$, or in other words, $Z^{-1} \in \hat{R}$.
\item  $R$ and $\hat{R}$ are Artin rings \cite{Atiyah_69}. In particular, every nonzero $A \in R$ is either a unit or a zero divisor, and the same property holds for any nonzero $Z \in \hat{R}$.
\end{itemize}

\begin{lemma}\label{killmanyblocks}
Let $A_1 \in R$ be invertible and $A_i \in R$. Define $Z_{m-1}, W \in \hat{R}$ as follows:
\[Z_{m-1}=\begin{bmatrix}
A_0 & A_1 & A_2 & \dots & A_{m-1}\\
0 & A_0 & A_1 & \dots & A_{m-2}\\
\vdots & \ddots & \ddots & \ddots & \vdots\\
\vdots & & 0 & A_0 & A_1\\
0 & 0 & \dots & 0 & A_0
\end{bmatrix}, \ W=\begin{bmatrix}
A_0 & I_n & 0 & \dots & 0\\
0 & A_0 & I_n & \dots & 0\\
\vdots & \ddots & \ddots & \ddots & \vdots\\
\vdots & & 0 & A_0 & I_n\\
0 & 0 & \dots & 0 & A_0
\end{bmatrix}.\]
 Then $W$ and $Z_{m-1}$ are similar over $R$, i.e., they are similar and the similarity matrix can be taken to have block elements in $R$.
\end{lemma}
\begin{proof}
Let $Z_1 = I_m \otimes A_0 + N_m \otimes A_1$. Clearly, $Z_1$ and $W= I_m \otimes A_0 + N_m \otimes I_n$ are similar over $R$ via $\diag(A_1,(A_1)^2,\dots,(A_1)^m)$. Hence, it suffices to show that $Z_{m-1}$ and $Z_1$ are similar over $R$.
We will explicitly exhibit an algorithm to construct a similarity matrix $X$ such that $Z_{m-1} X = X Z_1$. More specifically, $X$ can be found having the block structure

\[ X=\begin{bmatrix}
I_n & X_{1,2} & X_{1,3} & \dots & X_{1,m}\\
0 & I_n & X_{2,3} & \dots & X_{2,m}\\
\vdots & & \ddots & & \vdots\\
0 & \dots & & I_n & X_{m-1,m}\\
0 & 0 & \dots & 0 & I_n
\end{bmatrix},\]
 with $X_{i,j} \in R$. To this goal we will solve the equation $0=D:=Z_{m-1} X - X Z_1$ over $R$, i.e., viewing each block matrix as a matrix with scalar elements in the ring $R$. Coherently with this approach, within this proof we will from now on refer to block columns, block rows, block elements, etc., simply as columns, rows, elements, etc. 
 
There is no loss of generality in assuming $A_0=0$. Indeed, we may equivalently solve $(Z_{m-1} - I_m \otimes A_0) X = X (Z_1 - I_m \otimes A_0)$ because every matrix in $R^{m \times m}$ commutes with $I_m \otimes A_0$.  Let $U = N_m \otimes I_n \in \hat{R}$ and, for $k=1,\dots,m-1$, define $Z_k = \sum_{i=1}^k A_i U^i$. (Note that this definition also covers $Z_1,Z_{m-1}$.) The action of $U$ when acting from the left is an upwards row shift, while the action of $U$ from the right is a rightwards column shift. Since $R$ is a commutative ring we have\footnote{Here and below the notation $X(a:b,c:d)$, where $a \leq b, c \leq d$ are positive integers, indicates the submatrix of $X$ obtained by selecting the rows with indices $a,\dots,b$ and the columns with indices $c,\dots,d$. } $X Z_1=A_1 \begin{bmatrix}
0 & X(1:m-1,1:m-1)\\
0 & 0
\end{bmatrix}$. 
On the other hand, 
$\displaystyle Z_{m-1} X=\sum_{i=1}^{m-1} A_i \begin{bmatrix}
0 & X(i+1:m,i+1:m)\\
0 & 0
\end{bmatrix}$ (
the blocks have different sizes for different terms in the summation).
Therefore, $D_{ij}=0$ if $j-i<2$, while
$\displaystyle D_{ij}= A_{j-i} - A_1 X_{i,j-1} + \sum_{\ell=1}^{j-i-1} (A_\ell X_{i+\ell,j})$  if $j-i \geq 2$.
Thus, $D_{ij}$ is an affine function of $A_1$: $D_{ij}=D_{ij}^{(0)}+D_{ij}^{(1)}A_1$. In particular, when $j \geq i+2$, $D_{ij}^{(1)}=X_{i+1,j} - X_{i,j-1}$ and $D_{i,j}^{(0)}={A_{j-i}+}\sum_{\ell=2}^{j-i-1} (A_\ell X_{i+\ell,j})$. Both $X_{i+1,j}$ and $X_{i,j-1}$ belong to the $(j-i-1)$th superdiagonal\footnote{The first superdiagonal of a square matrix is the set of the elements directly above the diagonal. Recursively, the $k$th superdiagonal is defined as the set of the elements directly above the $(k-1)$th superdiagonal.} of $X$. Moreover, the unknown elements of $X$ appearing in $D^{(0)}_{ij}$ must belong to one of the first $j-i-2$ superdiagonals of $X$. 
We conclude that the following procedure can be applied to solve the problem:
\begin{enumerate}
\item set $X_{1,j}=0$ for all $2 \leq j \leq m$;
\item for fixed $j-i=2,\dots,m$ consider the equations $D_{ij}=0$ corresponding to the $(j-i)$th superdiagonal of $D$: These are $m-j+i$ equations in $m-j+i$ unknowns that can be solved via forward substitution for the elements in the $(j-i-1)$th superdiagonal of $X$  (this step requires the invertibility of $A_1$);
\item substitute into the next superdiagonals of $D$, and iterate until solved.
\end{enumerate}
\end{proof}

If one drops the assumption that $A_{1}$ is invertible, $Z_{m-1}$ and $W$ may be not even similar over $\mathbb{F}$. For example, if $n=m=3$, $A_0=N_3^2$, $A_1=2 N_3$, $A_2=-2 I_3$, then $W^2 \neq 0 = Z_{m-1}^2$.

\section{The bivariate Jordan problem}

\subsection{Generic case: Full solution}\label{sec:generic}

The problem of finding the JCF of $P(J_m(\lambda),J_n(\mu))$ becomes tractable when the vector of first Hasse derivatives is nonzero at $(\lambda,\mu)$, i.e., 
\begin{equation}\label{eq:nonzerograd}
\begin{bmatrix}
    p_x\\
    p_y
\end{bmatrix}: =  \begin{bmatrix}
    \fpartial_x p(\lambda,\mu)\\
    \fpartial_y p(\lambda,\mu)
\end{bmatrix} \neq 0.   
\end{equation}
 Theorem \ref{main} below solves the bivariate Jordan problem assuming \eqref{eq:nonzerograd} for all eigenvalue pairs.

\begin{theorem}\label{main}
For each Jordan block of $X$ of size $m$ with eigenvalue $\lambda$ and for each Jordan block of $Y$ of size $n$ with eigenvalue $\mu$, assume that $(\lambda,\mu)$ is such that \eqref{eq:nonzerograd} holds. Then, in the JCF of the matrix \eqref{pmat}:

\begin{itemize}

\item if $p_x \neq 0 \neq p_y$, {there are} $\min(m,n)$ Jordan blocks {associated with the eigenvalue $p(\lambda, \mu)$} of the following sizes: $m+n+1-2 \min(m,n), m+n+3-2 \min(m,n), \dots, m+n-1$;

{\item if $p_y=0$ but either $m=1$ or $p_x \neq 0$, let $r$ be the minimal value of $1 \leq k \leq n-1$ such that the evaluation at $(\lambda,\mu)$ of the $k$th Hasse derivative of \eqref{pscal} with respect to $y$ is nonzero (or, if all such derivatives are zero, $r=n$), and let $a, b$ be the unique integers such that $n=ar+b$ with $0 \leq b < r$. Then, there are $(r-b) \min(m,a) + b \min(m,a+1)$ Jordan blocks associated with the eigenvalue $p(\lambda, \mu)$ of the following sizes: $r-b$ copies of Jordan blocks of sizes $m+a+1-2 \min(m,a), m+a+3-2 \min(m,a), \dots, m+a-1$ and $b$ copies of Jordan blocks of sizes $m+a+2-2 \min(m,a+1), m+a+4-2 \min(m,a+1), \dots, m+a$;

\item if $p_x = 0$ but either $n=1$ or $p_y \neq 0$, let $r$ be the minimal value of $1 \leq k \leq m-1$ such that the evaluation at $(\lambda,\mu)$ of the $k$th Hasse derivative of \eqref{pscal} with respect to $x$ is nonzero (or, if all such derivatives are zero, $r=m$), and let $a, b$ be the unique integers such that $m=ar+b$ with $0 \leq b < r$. Then, there are $(r-b) \min(n,a) + b \min(n,a+1)$ Jordan blocks associated with the eigenvalue $p(\lambda, \mu)$ of the following sizes: $r-b$ copies of Jordan blocks of sizes $n+a+1-2 \min(n,a), n+a+3-2 \min(n,a), \dots, n+a-1$ and $b$ copies of Jordan blocks of sizes $n+a+2-2 \min(n,a+1), n+a+4-2 \min(n,a+1), \dots, n+a$.}

% \item if $p_y=0$ but either $m=1$ or $p_x \neq 0$, $(r-b) \min(m,a) + b \min(m,a+1)$ Jordan blocks of the following sizes: $r-b$ copies of Jordan blocks of sizes $m+a+1-2 \min(m,a), m+a+3-2 \min(m,a), \dots, m+a-1$ and $b$ copies of Jordan blocks of sizes $m+a+2-2 \min(m,a+1), m+a+4-2 \min(m,a+1), \dots, m+a$, where
% \begin{itemize}
% \item  if for some $1 \leq k \leq n-1$ the evaluation at $(\lambda,\mu)$ of the $k$th Hasse derivative of \eqref{pscal} with respect to $y$ is nonzero, then $r$ is the minimal value of $k$ satisfying such a property, otherwise $r=n$;
% \item $a$ and $b<r$ are such that $n=ar+b$.
% \end{itemize}

% \item if $p_x = 0$ but either $n=1$ or $p_y \neq 0$, $(r-b) \min(n,a) + b \min(n,a+1)$ Jordan blocks of the following sizes: $r-b$ copies of Jordan blocks of sizes $n+a+1-2 \min(n,a), n+a+3-2 \min(n,a), \dots, n+a-1$ and $b$ copies of Jordan blocks of sizes $n+a+2-2 \min(n,a+1), n+a+4-2 \min(n,a+1), \dots, n+a$, where
% \begin{itemize}
% \item if for some $1 \leq k \leq m-1$ the evaluation at $(\lambda,\mu)$ of the $k$th Hasse derivative of \eqref{pscal} with respect to $x$ is nonzero, then $r$ is the minimal value of $k$ satisfying such a property, otherwise $r=m$;
% \item $a$ and $b<r$ are such that $m=ar+b$;
%\end{itemize}
\end{itemize}
\end{theorem}
The first item in Theorem \ref{main}, i.e., the case $p_x \neq 0 \neq p_y$, is not a new result. It is indeed equivalent to \cite[Corollary 1]{Norman_97}.

\begin{remark}\label{ifmis1}
If $m=1$ then the value of $p_x$ is not important. This is explicitly stated in Theorem \ref{main} for the case $p_y=0$. If $p_y \neq 0$, note that if $m=1$ there are no positive integers $\leq m-1$. In the third item of Theorem \ref{main}, this means $r=1=m=a$ and $b=0$. Hence, when $m=1$ and $p_y \neq 0$, the first and third items predict the same Jordan structure. Analogous observations can be made if $n=1$, switching the roles of $x$ and $y$. Hence, Theorem \ref{main} can be applied if either $(\lambda,\mu)$ is not a critical point or $(m-1)(n-1)=0$.
\end{remark}

Let us illustrate Theorem \ref{main} with a couple of examples.

\begin{example}\label{ex:tbr}
Consider, as in \cite[Research Problem 3.11]{Higham_book}, the Fr\'{e}chet derivative of $f(W)=W^2$ for $W=J_2(0)=\begin{bmatrix}
0 & 1\\
0 & 0
\end{bmatrix}$. The associated bivariate polynomial is $p(x,y)=x+y$, therefore $p_x=p_y=1$. Since there is a single Jordan block of size $2$ and eigenvalue $0$, Theorem \ref{main} correctly predicts that the JCF of \eqref{eq:kron} is $J_3(0)+J_1(0)$.
\end{example}

 In fact, $p(x,y)=x+y$ as in Example \ref{ex:tbr} corresponds more generally to an important special case of Theorem \ref{main}: The JCF of a Kronecker sum. When specializing to this case, from Theorem \ref{main} we recover a classical result that we recall below.

\begin{theorem}[Jordan canonical form of a Kronecker sum]\label{thm:kronsum}\cite[Theorem 2]{Roth_34}
For each Jordan block of $X$ of size $m$ with eigenvalue $\lambda$ and for each Jordan block of $Y$ of size $n$ with eigenvalue $\mu$, {in} the JCF of the Kronecker sum $X \otimes I + I \otimes Y$ there are, associated with each eigenvalue $\lambda+\mu$, $\min(m,n)$ Jordan blocks of sizes $m+n+1-2 \min(m,n), m+n+3-2 \min(m,n), \dots, m+n-1$.
\end{theorem}

\begin{example}
Let $X=J_2(0) \oplus J_1(1)$ and $Y=
J_2(2) \oplus J_1(3)$, where $\oplus$ denotes, here and throughout the paper, the direct sum of two matrices. We study the JCF of $P(X,Y)$ with $p(x,y)=y - 2x + xy - y^2$. From the interaction of the Jordan block with eigenvalue $2$ in $Y$ and any Jordan block in $X$, we get $p_x =0 \neq p_y$; in the other cases we get $p_x \neq 0 \neq p_y$. Carefully applying the appropriate cases in Theorem \ref{main} to the four possible interactions of a Jordan block in $X$ with one in $Y$, we obtain that the JCF of $P(X,Y)=I_3 \otimes Y - 2 X \otimes I_3 + X \otimes Y - I_3 \otimes Y^2$ is
\[ J_2(-2) \oplus J_2(-2) \oplus J_2(-2) \oplus J_1(-5) \oplus J_2(-6).  \]
\end{example}

We now proceed to prove Theorem \ref{main}. By the discussion in Section \ref{sec:simp}, it suffices to study the JCF of $P(J_m(\lambda),J_n(\mu))$ under the assumption \eqref{eq:nonzerograd}.
\subsubsection{The Jordan canonical form of $P(J_m(\lambda),J_n(\mu))$ for $p_x \neq 0$, $p_y \neq 0$}\label{sec:nn}

As discussed, this case appeared as \cite[Corollary 1]{Norman_97}. To make the paper self-contained, we nevertheless give an independent proof. Since $p_x \neq 0$, we can apply Lemma \ref{killmanyblocks} and claim that it suffices to find the JCF of $I_m \otimes B + N_m \otimes I_n$, where $B$ is defined in \eqref{ourmatrix}. Then, since
$p_y \neq 0$, we can apply Lemma \ref{killmanyblocks} again to transform $B$ to a bidiagonal form via a certain similarity matrix $X$. Hence, we reduce to the problem of finding the JCF of 
\[ I_m \otimes A + N_m \otimes I_n = \begin{bmatrix}
A & I_n & 0 & \dots & 0\\
0 & A & I_n & \ddots  & \vdots\\
\vdots & \ddots & \ddots & \ddots & 0\\
\vdots & & \ddots & \ddots & I_n \\
0 & 0 & \dots & 0 & A\end{bmatrix}, \]
where $A$ is bidiagonal Toeplitz, and its diagonal (resp., first superdiagonal) contains the element $p(\lambda,\mu)$ (resp., 1).

The matrix above is a Kronecker sum, and it is known that the lengths of its Jordan chains are $m+n+1-2k$ for $k=1,2,\dots,\min(m,n)$. This is stated in \cite[Theorem 2]{Roth_34} (a crucial claim in \cite{Roth_34} is verified in \cite{Marcus_75}). 
An alternative approach is discussed in \cite{Brualdi_85}, where it is shown that this problem is related to the study of the $k$-paths of the Cartesian product of $\Gamma_m$ and $\Gamma_n$, where $\Gamma_k$ is the digraph having vertices $1,2,\dots,k$ and arcs $12,23,\dots,(k-1)k$. We refer the reader to \cite{Brualdi_85} for definitions and details. A slight modification to the proof of \cite[Lemma 4.3]{Brualdi_85} leads to the same conclusion as in \cite[Theorem 2]{Roth_34}. We also note that, substituting every occurrence of $I_n$ in the first block superdiagonal with $I_n+N_n$ via a block-diagonal similarity, we could equivalently have used the results concerning Kronecker products: see \cite[Lemma 4.3]{Brualdi_85}, and also \cite{Aitken_34, Littlewood_35, Marcus_75, Roth_34}. Therefore the first item in Theorem \ref{main} is proved.

\subsubsection{The Jordan canonical form of $P(J_m(\lambda),J_n(\mu))$ for $p_x \neq 0=p_y$}\label{sec:nz}

We now address the case where only one of the two first order Hasse  derivatives is zero. There is no loss of generality in assuming that the nonzero one is $p_x$, as otherwise we can just swap the roles of the triples $(x,\lambda,m)$ and $(y,\mu,n)$ via Lemma \ref{switch}. Following the same approach as above, we therefore aim to find the JCF of

\[ I_m \otimes B + N_m \otimes I_n = \begin{bmatrix}
B & I_n & 0 & \dots & 0\\
0 & B & I_n & \ddots  & \vdots\\
\vdots & \ddots & \ddots & \ddots & 0\\
\vdots & & \ddots & \ddots & I_n \\
0 & 0 & \dots & 0 & B\end{bmatrix}. \]
where $B$ is now an upper triangular Toeplitz matrix whose first superdiagonal is $0$. We now need a variation on the theme of Lemma \ref{killmanyblocks}.

\begin{lemma}\label{killnotasmanyblocks}
Let $\{A_i\}_{i=0}^{m-1} \subset R$. Define $Z_{m-1} \in \hat{R}$ by prescribing its first row to have elements $A_0, 0=A_1=\dots=A_{r-1},A_r,\dots,A_{m-1}$. That is, 
\[Z_{m-1}=\begin{bmatrix}
A_0 & 0 & \dots & 0 & A_r & A_{r+1}&\dots & A_{m-1}\\
0 & A_0 & 0 & \dots & 0 & A_r & \dots & A_{m-2}\\
\vdots & \ddots & \ddots & \ddots & & \ddots & \ddots & \vdots \\
0&\dots&0&A_0&0&\dots&&A_r\\
\vdots&0&\dots&0&A_0&0&\dots&0\\
\vdots&&0&\dots&0&A_0&\ddots&\vdots\\
\vdots&&&\ddots&&\ddots&\ddots&0\\
0 & \dots & \dots & \dots & 0 & \dots &0& A_0\\
\end{bmatrix}.\]
Suppose that $A_r \in R$, $r \leq m-1$, is a unit of $R$ (that is, it is invertible as a matrix over $\F$). Denote by $N_m$ the $m \times m$ Jordan block associated with the eigenvalue $0$. Then $Z_{m-1}$ and $I_m \otimes A_0 + N_m^r \otimes I_n$, are similar over $R$, i.e., they are similar and the similarity matrix can be taken having block elements in $R$.
\end{lemma}

\begin{proof}
We sketch the proof, that is constructive and similar to the one of Lemma \ref{killmanyblocks}. One  can find a similarity matrix $X \in R^{m \times m}$, upper triangular with all diagonal elements equal to $I_n$ and the first $r$ (block) rows equal to the rows of an identity matrix, such that $Z_{m-1} X = X(I_m \otimes A_0 + N_m^r \otimes A_r)$. It is then easy to find a (block) diagonal similarity matrix $D$ such that $D^{-1}(I_m \otimes A_0 + N_m^r \otimes A_r)D=(I_m \otimes A_0 + N_m^r \otimes I_n)$.
\end{proof}

Using Lemma \ref{killnotasmanyblocks} (in the scalar case $R=\F$) we see that it suffices to compute the JCF of the nilpotent matrix $N=I_m \otimes N_n^r + N_m \otimes I_n$. Let $q(x,y)=x+y$; then $N=Q(N_m,N_n^r)$. Note that $q_x = q_y = 1$ and therefore we can apply the results obtained in the previous section. The JCF of $N_m$ is $N_m$ itself, while the one of $N_n^r$ is provided by  Corollary \ref{jobloofpow}, which is a straightforward consequence of \cite[Theorem 1.36]{Higham_book}.

\begin{corollary}\label{jobloofpow}
Let $0 \leq b \leq r-1$ such that $n \equiv_r b$ and $n=ar+b$. Then, the JCF of the matrix $N_n^r$ has:
\begin{itemize}
\item $r-b$ Jordan blocks of size $a$, and
\item $b$ Jordan blocks of size $a+1$.
\end{itemize}
\end{corollary}
 
Hence, we conclude that there are $r-b$ copies of Jordan chains of length $m+a+1-2k$ for $k=1,2,\dots,\min(m,a)$ and $b$ copies of Jordan chains of length $m+a+2-2k$ for $k=1,2,\dots,\min(m,a+1)$. Finally, note that if $m=1$ there is only one block in \eqref{ourmatrix}, so that there is no need to assume that $p_x \neq 0$ to obtain the result.

\subsection{Degenerate case: General bounds}\label{sec:bounds}

When \eqref{eq:nonzerograd} is false, an exhaustive analysis of the bivariate Jordan problem is a nontrivial combinatorial problem. In \cite{Norman_97}, C. Norman fully solved the case $m,n \leq 3$, but this strategy cannot be pushed much further. Indeed, as $m$ and $n$ are unbounded, so is the number of subcases one would have to consider. Moreover, and contrary to the optimistic expectations that one might get from the solution of the generic case, the JCF does not only depend on the nonzero pattern of the second Hasse derivatives of $p$, but also on the values of the nonzero coefficients. For a counterexample, take $X=Y=N_3$, and let $p_1=-x^2+2xy+y^2$ and $p_2=x^2+xy+y^2$. Then, even though $p_1$ and $p_2$ have the same zero pattern of Hasse derivatives evaluated at $(0,0)$, one can check that  $p_1(X,Y)$ is similar to $N_3 \oplus N_2 \oplus N_2 \oplus N_1 \oplus N_1$, while $p_2(X,Y)$ is similar to $N_3 \oplus N_2 \oplus N_1 \oplus N_1 \oplus N_1 \oplus N_1$. We emphasize that only $p_2$ corresponds to a formal Fr\'{e}chet derivative; we shall study the Fr\'{e}chet-Jordan problem more thoroughly in Section \ref{sec:further}.

The previous observation is frightening, and we will not delve further into the epic undertaking of a complete solution to the bivariate Jordan problem. Instead, we give some general bounds on the size and number of the Jordan blocks in the degenerate case when \eqref{eq:nonzerograd} fails, based on the the concept of \emph{local degree}, that we define next. Given a bivariate polynomial $p(x,y)$ as in \eqref{pscal}, we define the \emph{local degree of $p$ at $(\lambda,\mu) \in \F^2$} as the smallest positive integer $d \geq 1$ such that  there exists a biindex $\alpha=(\beta,\gamma) \in \N^2$ with $|\alpha|=\beta+\gamma=d$ and $\fpartial_{z^\alpha} p(\lambda,\mu) = \fpartial_{x^\beta y^\gamma}p(\lambda,\mu) \neq 0$. Note that $d \in \N$ except if all the Hasse derivatives of $p$ are zero, or equivalently, except if $p$ is constant. As the bivariate Jordan problem is trivial when $p$ is a constant, we henceforth exclude this case.

Recall Remark \ref{ifmis1}: If either $m=1$ or $n=1$ (or both), Theorem \ref{main} can be applied regardless of the value of $p_x$ (if $m=1$) or regardless of the value of $p_y$ (if $n=1$). For this reason, from now on we also assume $m,n>1$. {Recall also that $\lceil \cdot \rceil$ and $\lfloor \cdot \rfloor$ denote the ceil and floor functions respectively.}

\begin{proposition}\label{zerodivisor2}
Suppose that $p(x,y)$ is not constant and let $d \geq 1$ be the local degree of $p(x,y)$ at $(\lambda,\mu)$. Write $P(J_m(\lambda),J_n(\mu))=p(\lambda,\mu) I_{mn}+Z_{m-1}$, where $Z_{m-1}$ is nilpotent. Then the size of the Jordan blocks of $P(J_m(\lambda),J_n(\mu))$ is bounded above by $\lceil \frac{m+n-1}{d} \rceil$.
\end{proposition}

\begin{proof}
We will show that $Z_{m-1}^{k}=0$ for $k \geq \lceil \frac{m+n-1}{d} \rceil$. Recall from Section \ref{sec:simp} that the commutative ring of $m \times m$ triangular Toeplitz matrices with elements in a commutative ring $R$ is isomorphic to the quotient ring $R[x]/\langle x^m \rangle$. Therefore, in order to analyze a power of $Z_{m-1}$ it suffices to consider the scalar polynomial $(A_0 + A_1 x + \dots + A_{m-1} x^{m-1})^k \bmod x^m$. Before quotienting, by the multinomial theorem we have
\[ (A_0 + A_1 x + \dots + A_{m-1} x^{m-1})^k = \sum_{|\alpha|=k} \frac{k!}{\alpha_0 ! \alpha_1! \cdots \alpha_{m-1}!} \prod_{i=0}^{m-1} A_i^{\alpha_i} x^{i \alpha_i},\]
where $\alpha=(\alpha_0,\dots,\alpha_{m-1})$ is a multiindex equipped with the usual multiindex $1$-norm {$|\alpha|=\sum_{i=0}^{m-1} \alpha_i$}.  Let us study the term in $x^h$ in the expansion above, $h \leq m-1$. By definition of $d$ and \eqref{ourmatrix}, and since $R \cong \F[y]/\langle y^n \rangle$, it holds for all $j=0,\dots,d-1$ that $A_j \in R$ corresponds to an element $[p_j] \in \langle [y]^{d-j}\rangle \subset \F[y]/\langle y^n \rangle$. Therefore, $\displaystyle \prod_{j=0}^{d-1} A_j^{\alpha_j}=0$ as long as $\displaystyle \sum_{j=0}^{d-1} (d-j) \alpha_j \geq n$. We now aim to minimize the value of $\sum_{j=0}^{d-1} (d-j) \alpha_j$ subject to $|\alpha|=k$ and $\sum_{i=0}^{m-1} i \alpha_i = h$. By the constraints we obtain \[ dk - h = \sum_{j=0}^{m-1} (d-j)\alpha_j \leq \sum_{j=0}^{d-1} (d-j) \alpha_j.\]

Hence, when $dk \geq m+n-1$, all the coefficients of $x^h$ are zero for all $h=0,\dots,m-1$.
\end{proof}
\vskip2pt
\begin{remark}
    Taking $d=1$ in Proposition \ref{zerodivisor2} yields \cite[Corollary 4.3]{Kressner_14}.
\end{remark}

 It is now convenient to introduce the following notation and terminology. Let $\{e_i\}_{i=1,\dots,m}$ denote the canonical basis of $\mathbb{F}^m$ and let $\{f_i\}_{i=1,\dots,n}$ be the canonical basis of $\mathbb{F}^n$. 
We define $V_0 = \{0\}$ 
and, for $1 \leq k\leq m+n-1$, $V_k = \mathrm{span}(\{e_i \otimes f_j  \ |  \ i+j \leq k+1   \})$. Note that $(V_k)_k$ is a \emph{filtration} of $\F^{mn}$, that is, a nested sequence of subspaces ($V_\kappa \subseteq V_k$ for $\kappa \leq k$) such that $V_k=\mathbb{F}^{mn}$ for $k \geq m+n-1$. Assuming with no loss of generality $m \leq n$, by an elementary enumeration argument we obtain the following formulae that describe the dimensions of the quotient spaces $V_k/V_{k-1}$.
\begin{equation}\label{eq:formulae}
  u_j:=\dim V_j/V_{j-1} = 
\min(j, m, m+n-j) \quad \text{for } j=1, \dots, m+n-1.
\end{equation}
 We say that $v \in \mathbb{F}^{mn}$ has \emph{$V$-degree} equal to $d$ if $v \in V_d$ but $v \not \in V_{d-1}$. Equivalently, $v \in V_d$ has $V$-degree $d$ if $[v] \neq [0]$ in $V_d/V_{d-1}$.

 \begin{lemma}\label{vdegree2}
Write $P(J_m(\lambda),J_n(\mu))=p(\lambda,\mu) I_{mn}+Z_{m-1}$, where $Z_{m-1}$ is nilpotent. For all $i=1,\dots,m$, $j=1,\dots,n$, the following properties hold:
\begin{itemize}
\item the $V$-degree of the $(n(i-1)+j)$th column of $P(J_m(\lambda),J_n(\mu))$ is at most $i+j-1$;
\item the $V$-degree of the $(n(i-1)+j)$th column of $Z_{m-1}$ is at most $\max(0,i+j-2)$;
\item if the local degree of $p$ at $(\lambda,\mu)$ is $d$,  then the $V$-degree of the $(n(i-1)+j)$th column of $Z_{m-1}$ is at most $\max(0,i+j-1-d)$.
\end{itemize}
\end{lemma}

\begin{proof}
By \eqref{ourmatrix}, we can write the following formula for the $(r,c)$ element of $P(J_m(\lambda),J_n(\mu))$ where $r=n(i_r-1)+j_r$ and $c=n(i_c-1)+j_c$. 

\begin{equation}\label{eq:thelements}
P(J_m(\lambda),J_n(\mu))_{rc} = \begin{cases}
  \fpartial_{x^{h}} \fpartial_{y^{k}} p(\lambda,\mu) \ &\mathrm{if} \ h:=i_c -i_r\geq 0 \ \mathrm{and} \ k:=j_c - j_r\geq 0;   \\
  0 \   &\mathrm{otherwise}.
\end{cases}    
\end{equation}
Hence, the $(n(i-1)+j)$th column of $P(J_m(\lambda),J_n(\mu))$ as $\displaystyle \sum_{h,k} q_{h k} N_m^h e_i \otimes N_n^k f_j$, with $q_{h k}=\left( \fpartial_{x^h} \fpartial_{y^k} \right) p$. If $q_{h k} \neq 0$ then the $(h,k)$ term in the summation has $V$-degree $i+j-h-k-1 \leq i+j-1$, yielding the first statement. The second (resp. third) statement is proven analogously, observing that in this case $q_{h k} = 0$ for $h=k=0$ (resp. $h+k \leq d-1$).
\end{proof}

Even if in this section we assume $m,n > 1$, Proposition \ref{prop:explicit_nullity} below is valid for any choice of $m$ and $n$.

\begin{proposition}\label{prop:explicit_nullity}
Let $d \ge 1$ be the local degree of $p(x,y)$ at $(\lambda,\mu)$. If $d \geq m+n-1$, then there are precisely $mn$ Jordan blocks in the JCF of $P(J_m(\lambda),J_n(\mu))$. If $d < m+n-1$, define $\delta = d - |n-m|$. Then, in the JCF of $P(J_m(\lambda),J_n(\mu))$, the number of Jordan blocks is bounded above by $\mathfrak{a} = \max(m,n) \min(m,n,d)$ and bounded  below by 
\[ \mathfrak{b} = \begin{cases}
d \min(m,n) \ &\mathrm{if} \ \delta \leq 0;\\
d \min(m,n) - \lfloor \frac{\delta^2}{4} \rfloor \ &\mathrm{if} \ \delta > 0.
\end{cases}\]
\end{proposition}

\begin{proof}
Write $P(J_m(\lambda),J_n(\mu))=p(\lambda,\mu) I_{mn} + Z_{m-1}$. If $d \geq m+n-1$, then  by \eqref{eq:thelements} every nonzero addend in $Z_{m-1}$ is a multiple of $N_m^a \otimes N_n^b$ with $a+b \geq m+n-1$. If $a \geq m$, $N_m^a=0$; and if $a < m$ then $b > n-1$ and $N_n^b=0$. Hence, $Z_{m-1}=0$.

Suppose now $d < m+n-1$. Recalling Lemma \ref{switch}, we may also assume with no loss of generality $m \leq n$. Since the number of Jordan blocks in the JCF of $P(J_m(\lambda),J_n(\mu))$ is the nullity of $Z_{m-1}$, our goal is to prove $\mathfrak{b} \leq \dim \ker Z_{m-1} \leq \mathfrak{a}$.

\begin{itemize}
\item \emph{Lower bound.} By Lemma \ref{vdegree2}, and recalling \eqref{eq:formulae}, for every $d \leq k \leq m+n-1$ there are $\displaystyle \sum_{j=1}^k u_j$ columns of $Z_{m-1}$ lying in $V_{k-d}$ (see also Figure \ref{fig1} for an illustration of this claim).
\begin{figure}
    \centering
   \[
\begin{bmatrix}
\mathbf{0} & \mathbf{0}  & \mathbf{1}  & \mathbf{2}  & \mathbf{0}  & \mathbf{1}  & \mathbf{2}  & \mathbf{3} & \mathbf{1} & \mathbf{2}  & \mathbf{3}  & \mathbf{4} \\
0 & 0 & \star & \star & 0 & \star & \star & \star & \star & \star & \star & \star\\
0 & 0 & 0 & \star & 0 & 0 & \star & \star & 0 & \star & \star & \star \\
0 & 0 & 0 & 0 & 0 & 0 & 0 & \star & 0 & 0 & \star & \star\\
0 & 0 & 0 & 0 & 0 & 0 & 0 & 0 & 0 & 0 & 0 & \star\\
0 & 0 & 0 & 0 & 0 & 0 & \star & \star & 0 & \star & \star & \star\\
0 & 0 & 0 & 0 & 0 & 0 & 0 & \star & 0 & 0 & \star & \star\\
0 & 0 & 0 & 0 & 0 & 0 & 0 & 0 & 0 & 0 & 0 & \star\\
0 & 0 & 0 & 0 & 0 & 0 & 0 & 0 & 0 & 0 & 0 & 0 \\
0 & 0 & 0 & 0 & 0 & 0 & 0 & 0 & 0 & 0 & \star & \star\\
0 & 0 & 0 & 0 & 0 & 0 & 0 & 0 & 0 & 0 & 0 & \star\\
0 & 0 & 0 & 0 & 0 & 0 & 0 & 0 & 0 & 0 & 0 & 0\\
0 & 0 & 0 & 0 & 0 & 0 & 0 & 0 & 0 & 0 & 0 & 0 \\
\end{bmatrix}.
\]
    \caption{Illustration of a claim within the proof of Proposition \ref{prop:explicit_nullity} for $n=4,m=3,d=2$. The symbol $\star$ denotes a generic, possibly nonzero, element of $\F$, and above each column we wrote in boldface the maximal possible $V$-degree.}
    \label{fig1}
\end{figure}
Hence, by the pigeonhole principle, $\displaystyle \dim \ker Z_{m-1} \geq \max_{d \leq k \leq m+n-1}  \ \sum_{j=1}^{k} (u_j - u_{j-d})$. We now have to distinguish two subcases, depending on whether $\delta$ is positive.
\begin{enumerate}
    \item Let $d \leq n-m$. In this case, $u_j < u_{j-d}$ precisely when $j > n$. The telescopic nature of the sum yields in turn the lower bound
    \[ \sum_{j=1}^{n} (u_j - u_{j-d}) = \sum_{j=n-d+1}^n u_j = dm.   \]
    \item Let $d>n-m$ so that $\delta > 0$. Then  $u_j \leq u_{j-d}$ precisely when $j \geq n + \left\lceil\frac{\delta}{2}\right\rceil$. Indeed,
\[m - \left\lceil \frac{\delta}{2} \right\rceil  = u_{n+\lceil \frac{\delta}{2} \rceil} \leq u_{n-d+\lceil \frac{\delta}{2} \rceil} = m - \left\lfloor \frac{\delta}{2} \right\rfloor\]
and
\[ m - \left\lceil \frac{\delta}{2} \right\rceil +1 = u_{n-1+\lceil  \frac{\delta}{2} \rceil} > u_{n-d{-1}+\lceil \frac{\delta}{2} \rceil} = m - \left\lfloor \frac{\delta}{2} \right\rfloor-1. \]
Therefore, in this case the lower bound is \[\displaystyle \sum_{j=m-\left\lfloor \frac{\delta}{2} \right\rfloor}^{n+\left\lceil \frac{\delta}{2} \right\rceil -1} u_j = dm -  \binom{\left\lceil \frac{\delta}{2} \right\rceil}{2} - \binom{1+\left\lfloor \frac{\delta}{2} \right\rfloor}{2}  = dm - \left\lfloor \frac{\delta^2}{4} \right\rfloor.  \]
We omit the tedious, but straightforward, verification of the arithmetic identity in the last step.
\end{enumerate}

    \item \emph{Upper bound.} If $d \geq \min(m,n)$, then $\mathfrak{a}=mn$ and thus the statement is trivial; hence, we assume $d < \min(m,n)$.  Since the number of Jordan blocks is the nullity of $Z_{m-1}$, in that case (still assuming $m \leq n$) we aim to prove $\mathrm{rank}(Z_{m-1}) \geq n(m-d)$. 
    
    By Lemma \ref{vdegree2}, for all $k$, $Z_{m-1}$ maps $V_k/V_{k-1}$ to $V_{k-d}$, where if $k<d$ we define $V_{k-d} =V_0= \{ 0 \}$. Let us permute rows and columns of $Z_{m-1}$, reordering them by increasing $V$-degree (for equal $V$-degrees, we keep the original order; in other words, after the permutation we order rows and columns by the grlex monomial order: $e_{i_1} \otimes f_{j_1} < e_{i_2} \otimes f_{j_1}$ when either $i_1+j_1 < i_2+j_2$ or $i_1+j_1=i_2+j_2$ and $i_1 < i_2$). We thus obtain a new matrix $T$ which is block upper triangular\footnote{Note that, generally, the blocks in a given superdiagonal of $T$ do not all have the same size.}. Furthermore, both the main block diagonal and the first $d-1$ block superdiagonals of $T$ are zero. Finally, $\displaystyle \mathrm{rank} (Z_{m-1}) = \mathrm{rank}(T) \geq \sum_{k=d+1}^{m+n-1} \mathrm{rank}(R_k)$ where $R_k$ are the blocks on the $(d+1)$-th block superdiagonal of $T$. In other words,  $\displaystyle R_k \in \F^{u_{k-d} \times u_k}$ is the banded (with band size $d+1$) Toeplitz rectangular matrix that represents the linear map  $V_{k}/V_{k-1} \rightarrow V_{k-d}/V_{k-d-1} $, obtained by appropriately restricting the action of $Z_{m-1}$. We know that at least one diagonal in the Toeplitz matrices $R_k$ is nonzero, and the ranks of $R_k$ are minimal when there is only one nonzero diagonal and it is in an extremal position. ({To see this, note that a lower bound for the rank of a banded Toeplitz matrix is given by the length of its extremal diagonals.)} In particular, since $m \leq n$, the worst-case scenario  occurs when  $\fpartial_{x^d} p(\lambda,\mu) \neq 0$ but the other derivatives of order $d$ are zero. Since we eventually have to sum these lower bounds over $k$,
\[ \mathrm{rank} (Z_{m-1}) \geq \mathrm{rank}( N_m^d \otimes I) = n \mathrm{rank}(N^d_m) = n (m-d).\] 
\end{itemize}
\end{proof}

We note that the proof of Proposition \ref{prop:explicit_nullity} can be adapted to obtain similar bounds on the nullity of the integer powers of $Z_{m-1}$. This can be done taking into account the following corollary.

\begin{corollary}
For all $i=1,\dots,m$, $j=1,\dots,n$, $k \geq 1$ we have the following properties:
\begin{itemize}
\item the $V$-degree of the $(n(i-1)+j)$th column of $Z_{m-1}^k$ is at most $\max(0,i+j-1-k)$;
\item if the local degree of $p(x,y)$ at $(\lambda,\mu)$ is $d$, then the $(n(i-1)+j)$th column of $Z_{m-1}^{k}$ has $V$-degree at most $\max(0,i+j-1-dk)$.
\end{itemize}
\end{corollary}

\begin{proof}
{We prove the second statement, which clearly implies the first}. The $(n(i-1)+j)$th column of $Z_{m-1}^k$ is $Z_{m-1}^k(e_i \otimes f_j)$. Note that by definition $(e_i \otimes f_j) \in V_{i+j-1}$. By Lemma \ref{vdegree2}, $Z_{m-1}  V_\ell \subseteq V_{\max(0,\ell-d)}$. Applying this statement $k$ times we obtain the result.
\end{proof}

\section{The Fr\'{e}chet-Jordan problem}\label{sec:further}
While the bivariate Jordan problem seems generally very hard, the Fr\'{e}chet-Jordan problem is more structured, and we focus on it in this section with the goal of going beyond the regime of \eqref{eq:nonzerograd}, that we already solved more generally.

We start with a simple observation. Let 
\begin{equation}\label{eq:chsbp}
   h_d(x,y):=\sum_{j=0}^d x^j y^{d-j}=\frac{x^{d+1}-y^{d+1}}{x-y} 
\end{equation}
denote the complete homogeneous symmetric bivariate polynomial of degree $d$.  Then, as observed in \cite[Theorem 5.1]{Kressner_14} the case of the formal Fr\'{e}chet derivative of the polynomial $f(w)$ corresponds to a very special subclass of bivariate polynomials:
\begin{equation}\label{eq:Bezout}
p(x,y)=\sum_{i=0}^{{L}} f_i h_i(x,y) = \frac{f(x)-f(y)}{x-y}. 
\end{equation}
Note that \eqref{eq:Bezout} is the B\'{e}zoutian of $f(w)$ and the constant polynomial $1$. The assumption that $p(x,y)$ has the form \eqref{eq:Bezout} imposes a rather rigid structure on the possible values of its formal Hasse derivatives, and thus on the JCF of $P(J_m(\lambda),J_n(\mu))$; in this section, we exploit this fact to solve the Fr\'{e}chet-Jordan problem. To explain the reason for this helpful rigidity, for a biindex $\alpha=(\beta,\gamma)$ define $z^\alpha:=x^{\beta}y^{\gamma}$ and $\fpartial_{z^\alpha}:=\fpartial_{x^{\beta}} \fpartial_{y^{\gamma}}$. Proposition \ref{xmenyxy} relates the derivatives of order $\ell-1$ of \eqref{eq:Bezout} to the mixed derivatives of order $\ell$, for all $\ell \geq 1$.

\begin{proposition}\label{xmenyxy}
Let $p(x,y)$ be of the form \eqref{eq:Bezout} for some $f(w) \in \F[w]$. Then, for all $\alpha \in \N^2$,
\[\left( \fpartial_{x z^{\alpha}}-\fpartial_{y z^{\alpha}} \right) p(x,y) = \left(x-y\right) \fpartial_{xyz^{\alpha}} p(x,y).\]
\end{proposition}
\begin{proof}
Multiplying \eqref{eq:Bezout} by $x-y$ yields $(x-y) p(x,y)=f(x)-f(y)$. The statement then follows by applying $\fpartial_{xyz^{\alpha}}$ to the latter equation.
\end{proof}
\vskip5pt
 Theorem \ref{mainfrechet}, which is ultimately based on Proposition \ref{xmenyxy}, solves the Fr\'{e}chet-Jordan problem over every algebraically closed field of characteristic zero.
\begin{theorem}\label{mainfrechet}
    Suppose that $p(x,y)$ has the form \eqref{eq:Bezout} and consider the matrix \eqref{eq:frechet}. Then, for each pair of Jordan blocks $\lambda I_m + N_m$ appearing in the JCF of $X$ and $\mu I_n + N_n$ appearing in the JCF of $Y$:
    \begin{enumerate}
        \item If $\lambda \neq \mu$, Theorem \ref{thm:distinct_ev} describes the corresponding Jordan blocks in the JCF of \eqref{eq:frechet};
        \item If $\lambda = \mu$, Theorem \ref{thm:equal_ev} describes the corresponding Jordan blocks in the JCF of \eqref{eq:frechet}.
    \end{enumerate}
\end{theorem}

The (rather involved) statements of Theorem \ref{thm:distinct_ev} and Theorem \ref{thm:equal_ev}, as well as their proofs, will be given in the next two subsections. Our approach differs depending on whether $\lambda$ and $\mu$ coincide or not. Note that, due to the form \eqref{eq:kron}, the case $\lambda=\mu$ is far from being a rare situation when studying the JCF of formal Fr\'{e}chet derivatives.

\subsection{Distinct eigenvalues}

Assume $\lambda \neq \mu$. Define the shifted matrix  $\widehat{P}(X,Y) = P(X,Y) - p(\lambda,\mu) I_{mn}$, corresponding to the scalar bivariate polynomial $\widehat{p}(x,y)=p(x,y)-p(\lambda,\mu)$.  Finally, define \begin{equation}\label{eq:phi}
   \varphi(w):=f(w)-w \ \frac{f(\lambda)-f(\mu)}{\lambda-\mu} \in \F[w], 
\end{equation}
where $f \in \F[w]$ is the polynomial appearing in the right hand side of \eqref{eq:Bezout}. Observe, in particular, that

\[ (\lambda-\mu) \begin{bmatrix}
    \fpartial_x p(\lambda,\mu)\\
    \fpartial_y p(\lambda,\mu)
\end{bmatrix} = \begin{bmatrix}
    \fpartial_w \varphi (\lambda)\\
    -\fpartial_w \varphi(\mu)
\end{bmatrix}. \]
Hence, \eqref{eq:nonzerograd} fails if and only if the Hasse derivative of $\varphi$ vanishes at both $\lambda$ and $\mu$.  Moreover, \eqref{eq:Bezout} yields

\begin{equation}\label{eq:Bezmatrix}
(J_m(\lambda) \otimes I_n - I_m \otimes J_n(\mu)) \cdot \widehat{P}(J_m(\lambda),J_n(\mu)) = \varphi(J_m(\lambda)) \otimes I_n - I_m \otimes \varphi(J_n(\mu)).   
\end{equation}

The left hand side of \eqref{eq:Bezmatrix} is the product of the invertible (because $\lambda\neq \mu$) matrix $J_m(\lambda) \otimes I_n - I_m \otimes J_n(\mu)$ and the nilpotent (by construction) matrix $\widehat{P}(J_m(\lambda),J_n(\mu))$; moreover, except for a shift of the eigenvalue, the latter matrix has the same JCF as $P(J_m(\lambda),J_n(\mu))$ In addition, by Remark \ref{rem:theycommute}, those matrices commute. These facts have interesting consequences.

\begin{lemma}\label{lem:commute}
Let $A \in \F^{k \times k}$ be invertible and let $B \in \F^{k \times k}$ be nilpotent. If $A$ and $B$ commute, then $AB$ and $B$ are similar.
\end{lemma}
\begin{proof}
Since $AB=BA$, we have $(AB)^k=A^k B^k$ for all $k \in \N$. In particular, $AB$ is nilpotent. The statement follows because, again for all $k$, $\rank(A^kB^k)=\rank B^k$.
\end{proof}

Theorem \ref{thm:distinct_ev} below characterizes the JCF of \eqref{eq:frechet}, in the case where $X$ and $Y$ are single Jordan blocks, with distinct eigenvalues. In its statement, we exclude the trivial case where $\varphi$ is {a constant}, i.e., $f$ is {an affine function} and $p$ is constant.

\begin{theorem}\label{thm:distinct_ev}
Suppose that $p(x,y)$ has the form \eqref{eq:Bezout}, and consider the matrix \eqref{eq:frechet}. Let $\lambda I_m + N_m$ and $\mu I_n + N_n$ be, respectively, a Jordan block in the JCF of $X$ and a Jordan block in the JCF of $Y$, such that $\lambda \neq \mu$. Define the positive integers $k,h \in \N$ as follows: 
\[ \fpartial_{w^i} \varphi(\lambda) = 0 \ \forall \ 1 \leq i < k, \quad \fpartial_{w^k} \varphi(\lambda) \neq 0,  \qquad \fpartial_{w^j} \varphi(\mu) = 0 \ \forall \ 1 \leq j < h, \quad \fpartial_{w^h} \varphi(\mu) \neq 0,   \]
where $\varphi \in \F[w]$ is the polynomial defined in \eqref{eq:phi} and we assume {$\varphi(w) \not\in \F$}. Moreover, define the integers $a,q,b,r$ as follows. If $k \geq m$ (resp. $h \geq n$), then $q=m$ and $a=0$ (resp. $r=n$ and $b=0$); otherwise, $m=a k + q$ (resp. $n=b h + r$) is the Euclidean division of $m$ by $k$ (resp. of $n$ by $h$). Finally, define the multisets $\{s_i\}_{i=1}^{\min(k,m)}$ and $\{ t_j \}_{j=1}^{\min(h,n)}$ where

\[ s_i = \begin{cases} a +1 \ &\mathrm{if} \  1 \leq i \leq q;\\
a \ &\mathrm{if}  \ q+1 \leq i \leq \min(k,m);\\
\end{cases} \qquad t_j = \begin{cases} b +1 \ &\mathrm{if} \  1 \leq j \leq r;\\
b \ &\mathrm{if}  \ r+1 \leq j \leq \min(h,n).\\
\end{cases}\]

Then, for each such pair of Jordan blocks of $X$ and $Y$, and for each of the $\min(k,m)\min(h,n)$ pairs $(i,j)$, with $1 \leq i \leq k$ and $1 \leq j \leq h$, in the JCF of \eqref{eq:frechet} there are  $\min(s_i,t_j)$ Jordan blocks of eigenvalue $\displaystyle \frac{f(\lambda)-f(\mu)}{\lambda-\mu}$ and sizes $s_i+t_j+1-2 \min(s_i,t_j), s_i+t_j+3-2 \min(s_i,t_j),\dots,s_i+t_j-1$.
\end{theorem}

\begin{proof}
    By \eqref{eq:Bezmatrix} and Lemma \ref{lem:commute}, except possibly for the eigenvalues, the JCF of $P(J_m(\lambda),J_n(\mu))$ is the same as that of $\varphi(J_m(\lambda)) \otimes I_n - I_m \otimes \varphi(J_n(\mu))$. The statement then follows by combining Theorem  \ref{thm:kronsum} and \cite[Theorem 1.36]{Higham_book}.
\end{proof}

Let us illustrate Theorem \ref{thm:distinct_ev} with a  couple of examples.

\begin{example}
Let us first explore an example not covered by Theorem \ref{main}. Consider the polynomial $f(w) = w^4 - 2w^2 \in \F[w]$ and the matrices $X = J_4(-1)$ and $Y = J_3(1)$. Then, $\varphi(w) = f(w)$, and it can be verified that $k=h=2$. The corresponding partitions are $\{s_1,s_2\}=\{2,2\}$ and $\{t_1,t_2\}=\{2,1\}$
Theorem \ref{thm:distinct_ev} predicts that 
the JCF of $P(X,Y)$ is \[ J_3(0) \oplus J_3(0) \oplus J_2(0) \oplus J_2(0) \oplus J_1(0) \oplus J_1(0) .\]
\end{example}

\begin{example}
We now consider a case where both Theorem \ref{main} and Theorem \ref{thm:distinct_ev} are applicable.
Let $f(w) = w^3 - w^2 \in \F[w]$, $X = J_4(0)$, and $Y = J_3(1)$; then $\varphi(w)=f(w)$, $k=2$, $h=1$, $\{s_1, s_2\} = \{2, 2\}$ and $\{t_1\} = \{3\}$. Note also that $k>1=h$ is equivalent to $p_x=0\neq p_y$. Both Theorem \ref{main} and Theorem \ref{thm:distinct_ev} predict that the JCF of $P(X,Y)$ is 
\[  J_4(0) \oplus J_4(0) \oplus J_2(0) \oplus J_2(0). \]
\end{example}

\subsection{Equal eigenvalues}
Let now $\lambda=\mu$. Similarly to the case $\lambda \neq \mu$, it is convenient to define a shifted version of $f$.
\begin{equation}\label{eq:phi2}
    \phi(w):=f(w) - w \ \fpartial_w f(\lambda).
\end{equation}
This guarantees that
\begin{equation}\label{eq:ecchenon}
  p(x,y) - p(\lambda,\lambda) = \frac{f(x)-f(y)}{x-y} - \fpartial_w f(\lambda) = \frac{\phi(x)-\phi(y)}{x-y}.     
\end{equation}

\begin{lemma}\label{lem:localdegfrechet}
    Let $p$ have the form \eqref{eq:Bezout} and let $\phi(w)$ be the polynomial defined in \eqref{eq:phi2}. Denote by $d$ the local degree of $p(x,y)$ at $(\lambda,\lambda)$ and by $\mathfrak{m}$ the multiplicity of $\lambda$ as a root of $\fpartial_w \phi(w)$. Then, $d=\mathfrak{m}$.
\end{lemma}
\begin{proof}
It is clear from \eqref{eq:phi2} that $\lambda$ is a root of $\fpartial_w \phi(w) = \fpartial_w f(w)- \fpartial_w f(\lambda)$. Hence, $\mathfrak{m} > 0$, coherently with the definition of local degree.   By Proposition \ref{xmenyxy}, and writing $z^\alpha:=x^\beta y^\gamma$ for every biindex $\alpha=(\beta,\gamma)$,  {$\displaystyle \fpartial_{w^{1+|\alpha|}}f(\lambda) = \fpartial_{z^\alpha} p(\lambda,\lambda)$.}
    Note now that $1 + |\alpha| \neq 0$ because $\mathrm{char}(\F) = 0$. It suffices to observe that, whenever $|\alpha|>0$, it holds \[ { \fpartial_{w^{|\alpha|}} \fpartial_w\phi(\lambda)=(1+|\alpha|) \fpartial_{w^{1+|\alpha|}} f(\lambda).} \]
\end{proof}

The proof of Lemma \ref{lem:localdegfrechet} also shows that, for Fr\'{e}chet derivatives and when $\lambda=\mu$, $p_x=p_y$ in the notation of Theorem \ref{main}. 

To state Lemma \ref{lem:formalroot}, denote by $\F[[w]]$ the ring of formal power series in $w$ with coefficients in $\F$, that is,
\[ \F[[w]] : = \left\{ \sum_{i=0} ^\infty c_i w^i, \ c_i \in \F \right\}. \]
\begin{lemma}\label{lem:formalroot}
    Let $\lambda \in \F$ and let $1 \leq k$ be an integer.  If $g(w) \in \F[w]$ satisfies $g(\lambda) \neq 0$, then there exists $r(w) \in \F[[w-\lambda]]$ such that the equation $[r(w)]^{k} = g(w)$ holds in $\F[[w-\lambda]]$.
\end{lemma}
\begin{proof}
    Since $g(w)$ is a polynomial, it is clear that it admits an expansion around $\lambda$, that is, $\displaystyle g(w) =  \sum_{i=0}^{\deg g} g_i (w-\lambda)^i$, with $g_0 \neq 0.$
   Let now $\eta \in \F$ be a $k$-th root of $g_0$;  the existence of $\eta$ is guaranteed because $\F$ is algebraically closed. We then see that it suffices to prove the statement when $g_0=1$, for if we can find $\hat{r}(w)$ such that $[\hat{r}(w)]^k = g(w)/g_0$ then we can take $r(w)=\eta \ \hat{r}(w)$ to guarantee $[r(w)]^k=g(w)$. Note now that the polynomial $y^k-1 \in \F[y]$ has nonzero discriminant, since it holds (using $\mathrm{char}(\F)=0$) $\mathrm{Res}(y^k-1,ky^{k-1}) = -(-k)^k \neq 0$. We can therefore invoke \cite[Proposition 6.1.9]{Stanley_book}, from which the statement follows immediately.
\end{proof}

We remark that $\mathrm{char}(\F)=0$ is crucial for Lemma \ref{lem:formalroot}, which can fail otherwise. For example, let $\mathrm{char}(\F)=p > 0$ and take $k=p$, $\lambda=0 \in \F$, and $g(w)=1+w \in \F[w]$. Yet, for every $r(w) \in \F[[w]]$,
\[ r(w)=\sum_{i=0}^\infty c_i w^i \Rightarrow [r(w)]^p = \sum_{i=0}^\infty c_i^p w^{ip} \neq 1 + w.  \]

\begin{theorem}\label{lem:attained}
    Suppose that $p(x,y)$ has the form \eqref{eq:Bezout}, that it is not constant, and that it has local degree $d \geq 1$. Then, {for all $\ell \in \N$, $[P(J_m(\lambda),J_n(\lambda)) - p(\lambda,\lambda)I_{mn}]^\ell$ is similar to $H_d(N_m,N_n)^\ell$}, where $h_d(x,y)$ is as in \eqref{eq:chsbp}.    
\end{theorem}

\begin{proof}
Let $d$ be the local degree of $p(x,y)$ at $(\lambda,\lambda)$ and let $\phi(w)$ be as in \eqref{eq:phi2}. We may assume $d < m+n-1$, for otherwise the statement reduces to $0=0$.

By Lemma \ref{lem:localdegfrechet}, and since $p(x,y)$ is not constant, we can then write $\phi(w)-\phi(\lambda)=(w-\lambda)^{1+d} g(w)$ for some polynomial $g \in \F[w]$ such that $g(\lambda) \neq 0$. Using Lemma \ref{lem:formalroot}, take $\hat{r}(w) \in \F[[w-\lambda]]$ to be a formal power series such that $[\hat{r}(w)]^{1+d}=g(w)$. Then, defining \[r(w):=(w-\lambda)\hat{r}(w),\] it holds \[\phi(w)-\phi(\lambda)=[r(w)]^{1+d} \quad \mathrm{and}  \quad \fpartial_w r(\lambda) = \sqrt[1+d]{g(\lambda)} \neq 0.\]  In addition, define $s(w)$ to be the polynomial obtained by truncating the series $r(w)$ in powers of $w-\lambda$ up to the term of order $m+n-1$. This implies, in particular, that $s(\lambda)=0$, $\fpartial_w s(\lambda) \neq 0$, $s(J_m(\lambda))=r(J_m(\lambda))$, and $s(J_n(\lambda))=r(J_n(\lambda))$. Finally, note the identity
\begin{equation}\label{eq:magia}
 { \left[\frac{\phi(x)-\phi(y)}{x-y}\right]^\ell =   \left[\frac{r(x)-r(y)}{x-y}\right]^\ell  \ \left[\frac{[r(x)]^{1+d}-[r(y)]^{1+d}}{r(x)-r(y)}\right]^\ell  .  }
\end{equation}
Recall the notation $\displaystyle h_d(x,y)=\frac{x^{d+1}-y^{d+1}}{x-y}$ and define $\displaystyle q(x,y):=\frac{s(x)-s(y)}{x-y}$. We conclude, by \eqref{eq:ecchenon} and \eqref{eq:magia}, that

\[ { \left[P(J_m(\lambda),J_n(\lambda)) - p(\lambda,\lambda) I_{mn} \right]^\ell = \left[Q(J_m(\lambda),J_n(\lambda)) \right]^\ell \cdot  \left[H_d(s(J_m(\lambda)),s(J_n(\lambda)) \right]^\ell  .     }   \]
{By construction, the only eigenvalue of $[Q(J_m(\lambda),J_n(\lambda))]^\ell$ is $(\fpartial_w s(\lambda))^\ell \neq 0$. Hence, by Lemma \ref{lem:commute}, $[P(J_m(\lambda),J_n(\lambda)) - p(\lambda,\lambda) I_{mn}]^\ell$ is similar to $[H_d(s(J_m(\lambda)),s(J_n(\lambda))]^\ell$.} In turn, since $\fpartial_w s(\lambda) \neq 0$ and $s(\lambda)=0$, for $k \in \{m,n\}$ $s(J_k(\lambda))$ is similar to $J_k(s(\lambda))=N_k$ by \cite[Theorem 1.36]{Higham_book}. At this point, we get the statement by the same argument as in the beginning of Section \ref{sec:simp}.
\end{proof}

By Theorem \ref{lem:attained}, the problem of computing the JCF of \eqref{eq:frechet} when $X,Y$ are two Jordan blocks of size $m$ and $n$ and with the same eigenvalues is therefore reduced to computing the ranks of $H_d(N_m,N_n)^\ell, \ell \in \N$, where $d$ is the local degree of \eqref{eq:Bezout}. Indeed, it is a classical result \cite{Weyr_85} that, for a nilpotent matrix $Z$, the number of Jordan blocks of size $s$ and eigenvalue $0$ is given by \[2 \dim \ker Z^s - \dim \ker Z^{s-1} - \dim \ker Z^{s+1}.\]

If $d\ell \geq m+n-1$, then clearly $H_d(N_m,N_n)^\ell=0$. Otherwise, following the same argument as in the proof of the upper bound in Proposition \ref{prop:explicit_nullity}, we can see that the nilpotent matrix $H_d(N_m,N_n)^\ell$ is permutation similar to a matrix $T$ that, when $d\ell < m+n-1$, has precisely  one nonzero block diagonal. The block elements of this nonzero block diagonal are the banded Toeplitz rectangular matrices $R_k \in \F^{(u_k-\ell d) \times u_k}$, for $d \ell + 1 \leq k \leq m+n-1$, and where $u_k$ are the integers defined by \eqref{eq:formulae}. Specifically, $R_k$ represents the action of $H_d^\ell(N_m,N_n)$ when restricted to a linear map $V_k/V_{k-1} \rightarrow V_{k-\ell d}/V_{k-\ell d-1}$. In particular, the nonzero entries of $R_k$ corresponds to the coefficients of the homogeneous polynomial $h_d(x,y)^\ell = \left( \sum_{i=0}^d x^i y^{d-i} \right)^\ell$. Hence, $R_k$ depends only on the quintuple $(k,m,n,d,\ell)$, and in particular, $R_k \in \Z^{(u_k-\ell d) \times u_k}$ is an integer banded Toeplitz matrix. (Note that, since $\F$ has characteristic zero, then necessarily $\Z \subsetneq \F$.) By some tedious, but straightforward, computations, one can see  that this description is equivalent to Definition \ref{def:Toeplitzmatrices} below.

\begin{definition}\label{def:Toeplitzmatrices}
    Fix the positive integers $m,n,d,\ell,k$, such that $m \leq n$, and $d \ell +1 \leq k \leq m+n-1$. Let the integers $u_k$ be defined as in \eqref{eq:formulae}. Furthermore, for $0 \leq i \leq \ell d$, define $\gamma_i$ by
    \[ \sum_{i=0}^{\ell d} \gamma_i z^i = (1 + z + \dots + z^d)^\ell,  \]
    and set $\gamma_i=0$ for every other value of $i$.
    
    Then, we define  the integer Toeplitz matrices $R_k \in \Z^{u_{k-\ell d} \times u_k}$ by setting $(R_k)_{ij}=\gamma_{j-i+c_k}$ where
   \[  c_k = \begin{cases}
        0 \ &\mathrm{if} \ k \leq n;\\
         k-n \ &\mathrm{if} \ n \leq k \leq n+\ell d;\\
          \ell d \ &\mathrm{if} \ k \geq n+\ell d.
    \end{cases}   \]
\end{definition}

If we denote the first column and row in $R_k$, respectively, by $e_{i_c} \otimes f_{j_c}$ and $e_{i_r} \otimes f_{j_r}$, then in Definition \ref{def:Toeplitzmatrices} the quantity $c_k = i_c - i_r$ describes the offset in the banded matrix $R_k$.  In Proposition \ref{prop:someproperties} we annotate some further properties of $R_k$.

\begin{proposition}\label{prop:someproperties}
Fix $k,m,n,d,\ell$ such that $m \leq n$ and $d \ell +1 \leq k \leq m+n-1$. Then, the banded Toeplitz matrix $R_k$ in Definition \ref{def:Toeplitzmatrices} satisfies the following properties:
\begin{enumerate}
    \item $0 \leq c_k \leq \ell d$;
    \item $u_{k-\ell d} \leq u_k + c_k \leq u_{k-\ell d} + \ell d$;
    \item The  top left element of $R_k$ is equal to $\gamma_{c_k} > 0$;
    \item The bottom right element of $R_k$ is equal to $\gamma_{u_k-u_{k-\ell d}+c_k} > 0$;
    \item Let $F_s \in \F^{s \times s}$ be the flip matrix (the matrix whose elements are $1$ on the main antidiagonal and $0$ elsewhere). Then, $F_{u_{k-\ell d}}R_k F_{u_k} = R^T_{\ell d + m + n -k}$, for all $d \ell +1 \leq k \leq m+n-1$. In particular, $R_k$ and $R_{\ell d+ m + n - k}$ have the same rank.
\end{enumerate}
\end{proposition}
\begin{proof}
    \begin{enumerate}
        \item Immediate by the definition of $c_k$.
        \item We distinguish three cases.
    \begin{itemize}
        \item Let first $k \leq n$. Then $u_{k-\ell d} \leq u_k$ because $m\leq n$ (see also the proof of the lower bound in Proposition \ref{prop:explicit_nullity}). By \eqref{eq:formulae} it is clear that $|u_{k+1}-u_k|\leq 1$, implying $u_{k} \leq u_{k-\ell d} + \ell d.$
        \item Let $n \leq k \leq n+\ell d$. Then, using \eqref{eq:formulae}, $u_k=m+n-k$. Moreover, since $k-\ell d \leq n$, either $u_{k-\ell d}=m \geq m-\ell d$ or, if $k-\ell d<m$, $u_{k-\ell d}=k-\ell d \geq m-\ell d$. Hence, 
        \[   u_{k-\ell d} \leq m = u_k + c_k  \leq u_{k-\ell d} + \ell d. \]
        \item Let $k \geq n+\ell d.$ Since $k \geq k - \ell d \geq n$, by \eqref{eq:formulae} we have $u_k=m+n-k$ and $u_{k-\ell d}=m+n-k+\ell d$, implying
        \[ u_{k-\ell d} = m+n-k+\ell d = u_k + \ell d \leq u_k + 2\ell d = u_{k-\ell d} + \ell d.    \]
    \end{itemize}
    \item By the first item, $0 \leq c_k \leq \ell d$ and thus $\gamma_{c_k}$ is a coefficient of the polynomial $(1 + z + \dots + z^d)^\ell$. By the multinomial theorem, it follows that $\gamma_{c_k}$ is a positive integer.
    \item By the second item $0 \leq u_{k}-u_{k-\ell d}+c_k \leq \ell d$. We can then follow the same argument as in the previous item.
    \item For simplicity of notation, set $h:=\ell d + m + n - k$ and $k'=k-\ell d$. It follows from \eqref{eq:formulae} that the sizes of $R_k$ and $R_{h}$ are compatible with the statement. We leave as an exercise to the reader to verify, from \eqref{eq:formulae} and Definition \ref{def:Toeplitzmatrices}, that $c_{h} + u_{k'} =  c_k + u_k$ holds for all $k$ in the range of definition of $R_k$.   Hence,
\[  (F_{u_{k'}}R_k F_{u_k})_{ij} = (R_k)_{u_{k'} + 1 - i,u_k+1-j} = \gamma_{i-j+u_k-u_{k'} + c_k} = \gamma_{i-j + c_{h} }  = (R_{h})_{ji}.   \]
    \end{enumerate}
\end{proof}

\begin{theorem}\label{thm:equal_ev}
Suppose that $p(x,y)$ has the form \eqref{eq:Bezout}, let $f \in \F[w]$ be the polynomial appearing in \eqref{eq:Bezout} and consider the matrix \eqref{eq:frechet} assuming without loss of generality that $m \leq n$. Let $d \geq 1$ be the multiplicity of $\lambda \in \F$ as a root of the polynomial ${\fpartial_w\phi} \in \F[w]$ as defined in Lemma \ref{lem:localdegfrechet}. For every positive integer $\ell$ and for every $d \ell +1 \leq k \leq m+n-1$, let $\rho(k,m,n,d,\ell)$ be the rank of the Toeplitz matrix  $R_k$ of Definition \ref{def:Toeplitzmatrices}. Let

\[ \mathfrak{N}_s := \begin{cases}
    0 \ &\mathrm{if} \ s=0;\\
    mn - \sum_{k=ds+1}^{m+n-1} \rho(k,m,n,d,s) \ &\mathrm{if} \ 0<sd < m+n-1;\\
    mn & \mathrm{otherwise}.
\end{cases}   \]
Then, for each pair of Jordan blocks in the JCFs of $X$  and $Y$, having sizes $m$ and $n$ respectively and sharing the same eigenvalue $\lambda$, in the JCF of \eqref{eq:frechet}:
\begin{itemize}
    \item if $d \geq m+n-1$, there are $mn$ Jordan blocks of size $1$ and eigenvalue $\fpartial_w f(\lambda)$;
    \item if $d < m+n-1$, there are precisely
\[ n_s : = 2\mathfrak{N}_s - \mathfrak{N}_{s-1} - \mathfrak{N}_{s+1}   \]
Jordan blocks having size $s \geq 1$ and eigenvalue $\fpartial_w f(\lambda)$.
\end{itemize}

\end{theorem}

\begin{proof}
    The statement follows by Theorem \ref{lem:attained}, the discussion above, and \cite{Weyr_85}.
\end{proof}

\begin{example}
    Let us revisit Example \ref{ex:tbr} to verify that Theorem \ref{thm:equal_ev} makes the same predictions as Theorem \ref{main}. Clearly $d=1$ and $m=n=2$. Let us now use Theorem \ref{thm:equal_ev}. For $s \geq 3$, $\mathfrak{N}_s=4$ is trivial. For $s=1$, we must sum the ranks of the Toeplitz matrices

    \[ R_2 = \begin{bmatrix}
        1 & 1
    \end{bmatrix}, \qquad  R_3=\begin{bmatrix}
        1\\
        1
    \end{bmatrix}   \]
to obtain $\mathfrak{N}_1=4-1-1=2.$ Similarly, when $s=2$ we subtract from $mn=4$ the rank of $R_3=2$, and hence $\mathfrak{N}_2=3$.

    We conclude that the sequence $(\mathfrak{N}_s)_{s \in \N}$ in Theorem \ref{thm:equal_ev} is
    \[ 0,2,3,4,4,4,\dots     \]
and we recover again the fact that, for this example, \eqref{eq:kron} has JCF
\[ J_3(0) \oplus J_1(0).\]
\end{example}

\begin{example}
    Consider the formal Fr\'{e}chet derivative of $f(w)=w^4-6w^2$, corresponding to $p(x,y)=h_3(x,y)-6 h_1(x,y)$.  Take as matrix argument $X=J_3(1)$ and $Y=J_2(1)$.  We see that $\phi(w)=f(w)+8w$ and hence $\partial_w \phi(w) =4(w+2)(w-1)^2$, so that {the local degree at $(1,1)$ is} $d=2$. Since $m>n$, we apply Lemma \ref{switch} and swap to $X=J_2(1)$, $Y=J_3(1)$, preserving the JCF. Note also that $\phi(w)-\phi(\lambda)=(w+3)(w-1)^3$, and therefore (as expected) the arguments in the proof of Theorem \ref{lem:attained} apply with $r(w)=(w-1)\sqrt[3]{w+3}$. A formal series expansion around $\lambda=1$ yields in turn $\displaystyle s(w)=\sqrt[3]{4}(w-1)+\frac{(w-1)^2}{6 \sqrt[3]{2}}$ (note that $\F$, being a closed field of characteristic zero, contains the algebraic closure of $\mathbb{Q}$). 
    
    We can then apply Theorem \ref{thm:equal_ev} with $m=2,n=3, d=2$.
    The only value $s$ for which the computation of $\mathfrak{N}_s$ is not trivial is $s=1$. In this case, we subtract from $mn=6$ the ranks of
    \[  R_3=\begin{bmatrix}
        1&1
    \end{bmatrix}, \qquad  R_4 = \begin{bmatrix}
        1\\
        1
    \end{bmatrix}.\]
Hence, 
\[ (\mathfrak{N}_k)_k=0,4,6,6,6 \dots   \]
and  therefore the sought JCF is
    \[ J_2(-8) \oplus  J_2(-8) \oplus J_1(-8) \oplus J_1(-8).   \]
\end{example}

\begin{example}\label{mortimer}
    Let $p(x,y)=h_4(x,y)$, corresponding to the formal Fr\'{e}chet derivative of $f(w)=w^5$. Take as matrix arguments $X=Y=J_4(0))$; the local degree at $(0,0)$ is $d=4$. We have $\mathfrak{N}_s=16$ for $s \geq 2$. When $s=1$, the quantity $mn-\mathfrak{N}_1$ is equal to the sum of the ranks of
    \[ R_5 = \begin{bmatrix}
        1 & 1 & 1
    \end{bmatrix}, \qquad R_6 =\begin{bmatrix}
        1&1\\
        1 &1
    \end{bmatrix}, \qquad R_7=\begin{bmatrix}
        1\\ 
        1\\
        1
    \end{bmatrix}.  \]
    This yields
    \[ (\mathfrak{N}_k)_k = 0, 13, 16, 16, 16 \dots   \]
    Thus, in this case the JCF of \eqref{eq:frechet} is
    \[ \bigoplus_{i=1}^{3} J_2(0) \oplus \bigoplus_{i=1}^{10} J_1(0). \] 

    Due to the nonmaximal rank of $R_6$, the lower bound of Proposition \ref{prop:explicit_nullity} is not attained: 
    \[ \dim \ker Z_{m-1} =\mathfrak{N}_1 =  13 > 12 = \mathfrak{b}. \]
\end{example}

\subsection{Ranks of $R_k$: A discussion and an open problem}

Theorem \ref{thm:equal_ev} gives an implicit solution to the Fr\'{e}chet-Jordan problem when $\lambda=\mu$, but it is in some sense less satisfactory than Theorem \ref{thm:distinct_ev} in that it does not provide a closed formula in terms of $n,m,d$. Instead, it reduces the problem to that of computing the ranks of certain integer Toeplitz matrices. We make several remarks from different angles.

    \begin{enumerate}
    \item \emph{Computational viewpoint}. Fix $m,n,d,\ell$. The matrix $H_d(N_m,N_m)^\ell$ has size $mn$. Instead, for $m \leq n$, there are at most $O(n)$ matrices $R_k$ having size at most $O(m)$. Moreover, by item 5 in Proposition \ref{prop:someproperties}, it suffices to compute the ranks of about half of them. Thus, computing and summing the ranks of $R_k$ is computationally more appealing than directly computing the rank of $H_d^\ell$.  
        \item  \emph{Linear algebraic viewpoint}. If (having fixed $m,n,\ell,d$) every matrix $R_k$ has maximal rank $\min(u_k,u_{k-\ell d})$, then the nullity of the matrix $H_d^\ell(N_m,N_n)$ achieves the lower bound $\mathfrak{b}$ in Proposition \ref{prop:explicit_nullity} (having replaced $d$ with $\ell d$ in the definition of $\delta$, and thus of $\mathfrak{b}$, in the statement of Proposition \ref{prop:explicit_nullity}). Indeed, in that case, $\mathrm{rank} (H_d(N_m,N_n)^\ell)$ is equal to
\[   \sum_{k=\ell d+ 1}^{m+n-1} \min(u_k,u_{k-\ell d}) = mn +\sum_{k=1}^{\ell d} u_k - \sum_{k=1}^{m+n-1} \max(0,u_k-u_{k -\ell d}) = mn -\mathfrak{b}(m,n,\ell d),\]
where, mutatis mutandis, the last step can be proved by the same argument as in the proof of the lower bound in Proposition \ref{prop:explicit_nullity}. This is a typical situation, i.e., it is what happens for most quadruples $(m,n,d,\ell)$. Nevertheless, as we already have seen in Example \ref{mortimer}, a loss of rank in some $R_k$ can happen, provided there are nontrivial solutions to the Diophantine linear equations prescribed by $R_k$. Below, we provide more insight via
Proposition \ref{prop:sufficient}, Example \ref{ex:trickycase}, and Example \ref{ex:trickiercase}.
\item \emph{Ring theoretical viewpoint}.  The nullity of $H_d(N_m,N_n)^\ell$ can also be interpreted in terms of bivariate polynomial ideals. Denote by $\mathcal{R}$ the quotient ring $\F[x,y]/\langle x^m,y^n \rangle$. By Remark \ref{rem:theycommute}, we can interpret $H_d(N_m,N_n)^\ell \in \mathcal{K}_{N_m,N_n}$ as the equivalence class $[h_d(x,y)^\ell] = [h_d(x,y)]^\ell \in \mathcal{R}$. Therefore, the ring isomorphism in Remark \ref{rem:theycommute} induces the vector space isomorphism
\[ \ker H_d(N_m,N_n)^\ell \cong \mathcal{I}_\ell:=\left\{ [p(x,y)] \in \mathcal{R}\ \mid \ p(x,y)h_d(x,y)^\ell \equiv 0 \bmod \langle x^m,y^n \rangle  \right\}.   \]
The set $\mathcal{I}_\ell$ is the \emph{annihilator ideal} of $[h_d(x,y)]^\ell$. We omit the straightforward verification that it is indeed an ideal of $\mathcal{R}$; note also that, although we have emphasized its dependence on $\ell$, it is clear that $\mathcal{I}_\ell$ also depends on $m,n,$ and $d$. If $\ell_0$ is the smallest integer such that $d \ell_0 \geq m+n-1$, we have the ascending chain of ideals \[ \langle [0] \rangle = \mathcal{I}_0 \subsetneq \mathcal{I}_1 \subseteq \dots \subseteq \mathcal{I}_{\ell_0} = \mathcal{R} = \mathcal{I}_{\ell_0+1} = \dots   \]
Given $(m,n,d)$, the Jordan form of $H_d(N_m,N_n)$ is thus completely determined by the sequence $\dim_\F \mathcal{I}_\ell = \dim_\F \ker H_d(N_m,N_n)^\ell$.
   \end{enumerate}

We now give, in Proposition \ref{prop:sufficient}, a sufficient condition for $R_k$ to lose rank. By item 5 in Proposition \ref{prop:someproperties}, we may without loss of generality assume that $R_k$ has at least as many rows as  columns. It is based on the insight that, for appropriate values of $m,n,d,\ell$, then $[x-y]^\ell$ belongs to the annihilator ideal of $[h_d(x,y)]^\ell$.

   \begin{proposition}\label{prop:sufficient}
       Fix $m \leq n, d, \ell$ and consider the matrices $R_k$ as in Definition \ref{def:Toeplitzmatrices}. Suppose with no loss of  generality that $\lceil \frac{m+n+\ell d}{2} \rceil \leq k \leq m+n-1  $, so that $u_{k-\ell d} \geq u_k$, i.e., $R_k$ does not have more columns than rows. Denote by $g_k:=(\ell d - c_k) \bmod d+1$ the remainder of $\ell d - c_k$ in the Euclidean division by $d+1$, so that $0 \leq g_k \leq d$. If (1) $u_k > \ell$ and (2) $ 0 < g_k < d - u_{k-\ell d} + 2$, then $R_k$ does not have full rank, i.e., $\mathrm{rank}(R_k) < u_k$.
   \end{proposition}
   \begin{proof}
       Let $0 \neq v \in \Z^{u_k}$ be the vector whose top elements are the coefficients of the polynomial $(x-y)^\ell$, ordered following the lexicographic order; this is possible by condition (1). By the Toeplitz structure of $R_k$, that encodes a convolution, $R_k v$ contains those coefficients of $q(x,y):=(x-y)^\ell h_d(x,y)^\ell = (x^{d+1}-y^{d+1})^\ell$ that are associated with the powers of $y$ between $\ell d - c_k$ and $\ell d-c_k + u_{k-\ell d}  - 1$. Write, by Euclidean division, $\ell d -c_k = q_k (d+1) + g_k$. By condition (2),
       \[ q_k ( d+ 1) < \ell d - c_k \leq \ell d - c_k + u_{k-\ell d} - 1 < q_k(d+1) + d + 1 {=} (1+q_k)(d+1), \]
       ensuring $R_k v = 0$ because the only nonzero coefficients of $q(x,y)$ are associated with the powers of $y$ that are exact multiples of $d+1$.
   \end{proof}

    \begin{example}\label{ex:trickycase}
        We illustrate the concrete application of Proposition \ref{prop:sufficient} for $d=3$ and $\ell=2$. Then, the $\gamma_i$ are defined as the coefficients of
        \[(1+z+z^2+z^3)^2 = 1 + 2 z + 3 z^2 + 4 z^3 + 3 z^4+ 2 z^5+ z^6.\]
        Take now $m=4,n=8,k=9$. Then, $\ell d +1 = 7 \leq k \leq 11=m+n-1$ and $c_9=1$. Moreover 
\[ u_9 = u_3=3, \qquad R_9 = \begin{bmatrix}
    \gamma_1 & \gamma_2 & \gamma_3\\
    \gamma_0 & \gamma_1 & \gamma_2\\
    \gamma_{-1} & \gamma_0 & \gamma_1
\end{bmatrix} = \begin{bmatrix}
    2&3&4\\
    1&2&3\\
    0&1&2
\end{bmatrix} \Rightarrow \rho(9,4,8,3,2)=2 < 3.\]
Note that the null vector is $\displaystyle \begin{bmatrix}
    1\\
    -2\\
    1
\end{bmatrix}$, which encodes the coefficients of $(x-y)^2$, and $(x-y)^2 \in \mathcal{I}_2$ for the chosen values of $m,n,d,\ell$.
    \end{example}

The idea of Proposition \ref{prop:sufficient} can be generalized by taking, more generally, polynomials in the ideal $\langle (x-y)^\ell \rangle$. Example \ref{ex:trickiercase} shows one such construction at work.

    \begin{example}\label{ex:trickiercase}
        We illustrate one instance of a rank drop that is not predicted by Proposition \ref{prop:sufficient}. Let $\ell=1$ so that $\gamma_i \in \{0,1\}$. Take also $d=2,m=n=6,$ and $k=7$. Then, $u_7=u_5=5$, $c_7=1$ and
        \[  R_7 = \begin{bmatrix}
            1 & 1 &0&0&0\\
            1&1&1&0&0\\
            0&1&1&1&0\\
            0&0&1&1&1\\
            0&0&0&1&1
        \end{bmatrix} \Rightarrow \begin{bmatrix}
            1\\
            -1\\
            0\\
            1\\
            -1
        \end{bmatrix} \in \ker R_7.  \]
        From the viewpoint of ring theory, the vector in the null space encodes the polynomial $xy(x^3+y^3)(y-x)$, whose equivalence class belongs to $\mathcal{I}_1$ for the chosen values of the variables, since \[ h_2(x,y)(x-y)(x^3+y^3)=x^6-y^6 \equiv 0 \bmod \langle x^6,y^6 \rangle .\]
    \end{example}

While these examples and partial results suggest that our insight of constructing polynomials in the annihilator ideals of $[h_d(x,y)]^\ell$ is a powerful idea, a complete characterization of the ranks of $R_k$ involves the solution of the following combinatorial-algebraic problem: Given $d,\ell,m,n$ construct all the possible homogeneous polynomials $v(x,y)$ such that $v(x,y)h_d(x,y)^\ell$ has sufficiently many consecutive zero coefficients (where the meaning of ``sufficiently many" is that arising from the previous discussion). We find that a full analysis is beyond the scope of the current manuscript, and hence we  conclude with the open Problem \ref{prob:beauh}, whose solution would allow us to refine the statement of Theorem \ref{thm:equal_ev}.

\begin{problem}\label{prob:beauh}
    Fully characterize the quadruples $(m,n,d,\ell)$ for which some $R_k$ in Definition \ref{def:Toeplitzmatrices} is rank deficient, and quantify exactly the amount of these rank drops.
\end{problem}

Problem \ref{prob:beauh} is related to certain deep modern results in commutative algebra; namely, it can be seen as a generalization (from linear to non-linear polynomials in the graded Artin ring $\mathcal{R}=\F[x,y]/\langle x^m,y^n \rangle$) of the problem of determining the membership, or lack thereof, of $[h_d(x,y)^\ell] \in \mathcal{R}$ to the \emph{non-Lefschetz locus} of $\mathcal{R}$ \cite[Definition 2.3]{Boij_18}.

\section{Conclusions and future outlook}\label{sec:conclusions}

We studied two nontrivial problems of combinatorial nature, the Fr\'{e}chet-Jordan problem and the more general bivariate Jordan problem. We solved the Fr\'{e}chet-Jordan problem by reducing it to computing either the Euclidean division of two numbers or to computing ranks of certain integer Toeplitz matrices. For the bivariate Jordan problem, we gave some partial results, complementing the existing literature. A full solution of the bivariate Jordan problem remains an open problem, and indeed not an easy one in the author's opinion. Fully characterizing the ranks of the integer Toeplitz matrices that appear in our solution of the Fr\'{e}chet-Jordan problem, or extending our analysis to fields of prime characteristic, are also possible directions for future research.

\section{Acknowledgements}

The late Nick Higham suggested the Fr\'{e}chet-Jordan problem to me when I was a postdoc in his group. Nick himself, Bruno Iannazzo, and Samuel Relton (as well as Giovanni Barbarino, a few years later) shared comments on an earlier version of this document, that did not go much beyond Theorem \ref{main}. This earlier version was in fact submitted, and an anonymous reviewer pointed me to \cite{Kressner_14}, which in turn cites \cite{Norman_97}. I retracted that version, and decided to resubmit some time later, after having obtained more results.


\begin{thebibliography}{3}

\bibitem{Aitken_34} { A. C. Aitken}, The normal form of compound and induced matrices. \emph{Proc. London Math. Soc.} 38:354--376 (1934).  

{\bibitem{Atiyah_69} {M. F. Atiyah and I. G. Macdonald}, \emph{Introduction to Commutative Algebra}, Addison-Wesley, 1969.}

\bibitem{Berger_77} { M.~~S.~Berger}, {\em Nonlinearity and functional analysis: Lectures on nonlinear problems in mathematical analysis}, Academic Press, New York City, NY, USA, 1977.

\bibitem{Boij_18} M. Boij, J. Migliore, R.M. Mir\'{o}-Roig and U. Nagel, The non-Lefschetz locus, \emph{J. Algebra} 505:288--320 (2018).

\bibitem{Brualdi_85} { R. A. Brualdi}, Combinatorial verification of the elementary divisors of tensor products. \emph{Linear Algebra Appl.} 71:31--47 (1985).

\bibitem{GS72} { I. Gohberg and A. Semencul}, On the inversion of finite Toeplitz matrices and their continuous analogs. \emph{Mat. Issled.} 2:201--233 (1972).

\bibitem{Higham_book} { N.~J.~Higham}, {\em Functions of Matrices: Theory and Computation}, SIAM, Philadelphia, PA, USA, 2008.

\bibitem{HJunder} { R.~A.~Horn and C.~R.~Johnson}, {\em Matrix Analysis}, Cambridge University Press, 2013 (Second edition).

\bibitem{HJ} { R.~A.~Horn and C.~R.~Johnson}, {\em Topics in Matrix Analysis}, Cambridge University Press, 1991.


\bibitem{Littlewood_35} { D.~E.~Littlewood}, On induced and compound matrices. \emph{Proc. London Math. Soc.} 39(2):370--381 (1935). 

\bibitem{Kressner_14} {D. Kressner}, Bivariate matrix functions. \emph{Oper. Matrices} 8(2):449--466 (2014).

\bibitem{Marcus_75} { M. Marcus and H. Robinson}, Elementary divisors of tensor products. \emph{Commun. ACM} 18:36--39 (1975).

\bibitem{Norman_97} {C. W. Norman}, On the Jordan form of a family of linear mappings. \emph{Linear Algebra Appl.}, 257:225–-241 (1997).

\bibitem{Roth_34} { W. E. Roth}, On direct product matrices. \emph{Bull. Amer. Math. Soc.} 40:461--468 (1934).

\bibitem{Stanley_book} {R. P. Stanley}, \emph{Enumerative Combinatorics}, Volume 2, Cambridge University Press, 1999.

\bibitem{Weyr_85} {E. Weyr}, R\'{e}partition des matrices en esp\`{e}ces et formation de toutes les esp\`{e}ces.  \emph{C. R. Acad. Sci. Paris} 100:966--969 (1885).

\end{thebibliography}
\end{document}